\documentclass{article}

\usepackage{amsmath}
\usepackage[francais]{babel}
\usepackage{amssymb}
\usepackage[applemac]{inputenc}
\usepackage{mathrsfs}
\newcommand{\sm}{\raisebox{2pt}{~\rule{6pt}{1.4pt}~}}

\def\F{{\mathbb F}}

\def\N{{\mathbb N}}
\def\Z{{\mathbb Z}}
\def\Sq{{\rm Sq}}

\def\R{{\mathbb R}}
\def\rg{{\rm rang}}

\def\deg{\mathop{\rm deg}}
\def\<{\langle}
\def\>{\rangle}
\def\q_{\underline{q}}

\def\Tors{\mathop{\rm Tors}}
\def\Corps{\mathop{\rm Corps}}
\def\Ab{\mathop{\rm Ab}}
\def\Et{\mathop{\rm Et}}

\def\GL{\mathop{\rm GL}}

\def\Gal{\mathop{\rm Gal}}

\def\Sym{\mathop{\rm Sym}}

\def\Res{\mathop{\rm Res}}
\def\Cor{\mathop{\rm Cor}}
\def\Inv{\mathop{\rm Inv}}
\def\gal{\mathop{\rm gal}}
\def\gr{\mathop{\rm gr}}
\def\fil{\mathop{\rm fil}}
\def\Mil{\mathop{\rm Mil}}

\def\Tr{\mathop{\rm Tr}}
\def\cst{\mathop{\rm cst}}
\def\norm{\mathop{\rm norm}}

\def\Spec{\mathop{\rm Spec}}
\def\la{\lambda}
\def\e{\varepsilon}
\def\n{\noindent}

\def\q{\underline{q}}

\def\cc{\mathscr{C}}

\begin{document}

   \medskip
     
  \centerline{   {\bf  Invariants cohomologiques mod 2 et invariants de Witt des groupes altern\'es}} 
  
  \medskip

  \centerline{Jean-Pierre Serre}

   \medskip 
  Nous nous proposons d'\' etendre aux groupes altern\' es $A_n$ les r\' esultats obtenus dans [Se 03] pour les groupes sym\' etriques $S_n$, \`a la fois pour les invariants cohomologiques mod 2 et pour les invariants de Witt. 
  
  \smallskip
   Le $\S$1 rappelle des d\' efinitions et r\' esultats standards sur ces invariants.  Le $\S$2 donne les relations qui se d\' eduisent des (ex)-conjectures de Milnor. 
   
   \smallskip
   
   Le $\S$3 \' enonce les principaux r\' esultats (th\' eor\` emes 3.4.1 et 3.4.2); ils sont d\' emontr\' es aux $\S$$\S$ 4, 5. La m\' ethode de d\' emonstration est analogue \`a celle utilis\' ee pour $S_n$ et plus g\' en\' eralement pour les groupes de Coxeter finis (cf. [Se 03], [Se 18], [GH 23]): un ``th\' eor\` eme de d\' etection'' ({\it splitting principle}) permet de remplacer le groupe \' etudi\' e (ici $A_n$) par un sous-groupe plus simple : un $2$-sous-groupe ab\' elien \' el\' ementaire.  On en d\' eduit que les invariants cohomologiques et les invariants de Witt de $A_n$ sont d\' etect\' es par les alg\`ebres \' etales de dimension $n$ qui sont produits d'alg\`ebres de dimension $1$ et d'alg\`ebres biquadratiques (th. 5.5.1); rappelons que, dans le cas de $S_n$, on a le m\^ eme \' enonc\' e, avec
   ``biquadratiques'' remplac\' e par ``quadratiques'', cf. [Se 03], $\S$$\S$ 25, 29. Cela entra\^ine que le module des invariants est un module libre de rang $1 + [\frac{n}{4}]$ sur l'alg\` ebre de cohomologie de $k$ ( dans le cas des invariants cohomologiques), et sur
   l'anneau de Witt de $k$ (dans le cas des invariants de Witt). Noter que $1 + [\frac{n}{4}]$ est \' egal au nombre des classes de conjugaison d'involutions de $A_n$;
 la situation est analogue \`a celle des groupes de Weyl, cf. [Se 18].
   
   \smallskip
   
   Le $\S$6 compl\`ete les pr\' ec\' edents en pr\' ecisant les propri\' et\' es des formes traces associ\' ees \`a $A_n$.

     Un Appendice d\' emontre une formule utilis\' ee au $\S$1.1.

    \smallskip
    \n {\it Convention.} Tous les corps consid\' er\' es sont suppos\' es de caract\' eristique $\neq 2$.

     \newpage   
     
  \n    {\sc $\S$1. Rappels et compl\' ements}  

\smallskip

      Dans ce qui suit, $k$ est un corps de caract\' eristique $\neq 2$ et $G$  est un groupe fini.

            \medskip 
 \n  {\bf 1.1. Cohomologie mod 2 de $k$}.   
    
    \smallskip
     Soit $k_s$ une clôture s\' eparable de
      $k$, et soit $\Gamma_k = \Gal(k_s/k)$. Le groupe $\Gamma_k$ a une structure naturelle de groupe topologique profini. Dans la suite, les homomorphismes entre groupes profinis seront tacitement suppos\' es  continus.
      
      On note $H(k)$ l'alg\`ebre de cohomologie de $k$ mod 2, autrement dit :
      $$ H(k) =\oplus_{r \geqslant 0} H^r(k) = \oplus_{r \geqslant 0} H^r(\Gamma_k,\F_2).$$
    \quad \ \   On a $H^0(k)= \F_2$ et $H^1(k) \simeq k^\times/k^{\times 2}$. Si $x \in k^{\times}$, 
      l'\' el\' ement correspondant de $H^1(k)$ est not\' e  $(x)$, ou bien $(x)_k $ si l'on d\' esire pr\' eciser $k$. 
      
    \smallskip  Dans le cas $x=-1$,
    on \' ecrit $e_k$ \`a la place de $(-1)_k$. On a :
      
      \smallskip      
      (1.1.1) \quad $z\!\cdot\!(-z)=0$, i.e. $z^2=e_kz$, pour tout $z\in H^1(k)$.
      
      \smallskip
     \n   Les conjectures de Milnor ([Mi 70]), d\' emontr\' ees par Voevodsky ([Vo 03.I], [Vo 03.II], [OVV 07],  voir aussi  [Me 10], [Mo 05]) entraînent que l'alg\`ebre $H(k)$ est engendr\' ee par  ses \' el\' ements de degr\' e $1$. D'apr\`es (1.1.1), cela implique que, si $ z \in H^d(k)$, avec $ \ d \geqslant 1$, on a:

     \smallskip
        
         (1.1.2) \quad $z^m =  e_k^{dm-d}z$ pour tout $m \geqslant 1$.
         
         \smallskip \n  D' où      $y^2z=z^2y$ pour tout $y\in H^d(k)$, i.e.:
         
     \smallskip (1.1.3) \quad     $(1+y)(1+z) = (1+y+z)(1+yz)$.
     
     \smallskip \n Si $t\in H^d(k)$, on d\' eduit de (1.1.3) :
     
     \smallskip  (1.1.4) \quad $(1+y)(1+z)(1+t) = (1+y+z+t)(1+yz + zt +tx)$.

     \smallskip \n La formule (1.1.4) s'\' ecrit aussi  $(1+y)(1+z)(1+t) = (1+a_1)(1+a_2)$,   où  $a_1, a_2$ sont les fonctions sym\' etriques \' el\' ementaires 
   de degr\' e 1 et 2 de $y,z,t$.  
   
   \smallskip
   
   \n Plus g\' en\' eralement, soient
   $y_1,...,y_n\in H^d(k)$, et soient  $s_m \in H^{md}(k)$
  leurs fonctions sym\' etriques \' el\' ementaires ($m=1,...$). On a  $ 1+\sum_m s_m = \prod_i(1+y_i).$ En appliquant de façon r\' ep\' et\' ee la formule (1.1.3) on obtient  apr\`es quelques calculs (voir les d\' etails dans l'Appendice \`a la fin du texte, $\S$A.4):

\medskip
   
   (1.1.5) \quad $\prod_i(1+y_i) = \prod_j(1 + s_{2^j}).$
   
    \smallskip
  \n Autrement dit, si un entier  $m \geqslant 0$ a pour \' ecriture dyadique
   $m = \sum_{\alpha \in A} 2^\alpha$, avec $A \subset \N$, on a :
   
    \smallskip
   (1.1.6) \quad $s_m = \prod_{\alpha \in A} s_{2^\alpha}.$

 \smallskip   
 \n   Soit $m' = \sum_{\alpha \in A'}$ un autre entier $\geqslant 0$. D\' efinissons l'entier  $m\bullet m' $ ({\it somme diminu\' ee} de $m$ et $m'$) par la formule :
   
   \smallskip
 (1.1.7) \quad $m\bullet m'  =  \sum_{\alpha \in A \cup A'} 2^\alpha =m+m'-\sum_{\alpha \in A \cap A'} 2^\alpha. $
 
 \medskip\n Les formules (1.1.2) et (1.1.6) entraînent :
   
  \smallskip  
   (1.1.8)\quad $ s_ms_{m'} = e_k^{d.|A \cap A'|} s_{m\bullet m' }.$

  \smallskip \n Par exemple, si $m =3, m' = 6$, on a $m \bullet m' =7$, d'où $s_3s_6 = e_k^{2d}s_7 =e_k^{2d}s_1s_2s_4.$
  
  \medskip \n {\it Carr\' es de Steenrod}.
  
  \smallskip Les carr\' es de Steenrod op\`erent sur $H(k)$  (comme sur
  toute alg\`ebre de cohomologie mod 2 d'un groupe profini, ou d'un groupe discret). Le fait que $H(k)$ soit engendr\' e par des \' el\' ements de degr\' e 1 entraîne la formule suivante, valable pour tout $z \in H^d(k)$ :
          
        \smallskip         (1.1.9) \quad  $\Sq(z) = (1+e_k)^dz$, 
        
     \smallskip   \n   où  $\Sq(z)$ est le carr\' e de Steenrod ``total" de $z$, 
 i.e.  $\sum_{0\leqslant i \leqslant d}\Sq^i(z)$. En particulier, on a $\Sq^1(z) = 0$ si $d$ est pair et $\Sq^1(z)=e_kz$ si $d$ est impair.

      \medskip
      
      \n {\bf 1.2. Anneau de Witt de $k$}.
      
          \smallskip
      
        Soit $\widehat{W}(k)$ l'anneau de Grothendieck-Witt de $k$, autrement dit le
        groupe de Grothendieck des classes d'isomorphisme de $k$-formes quadratiques;
        il a une structure naturelle de $\la$-anneau commutatif (produit tensoriel, puissances ext\' erieures). Il est engendr\' e par les classes $\<\alpha\>$ des formes quadratiques $x \mapsto \alpha x^2$ de rang 1. La classe d'isomorphisme d'une forme quadratique $\alpha_1X_1^2 + \dots + \alpha_nX^2_n$ est not\' ee $\<\alpha_1,\dots,\alpha_n\>$.
        
        Le quotient de  $\widehat{W}(k)$ par l'id\' eal engendr\' e par la        
       forme hyperbolique $\< 1,-1\>$ est l'anneau de Witt $W(k)$. On a un carr\' e cart\' esien: 
       
       \smallskip
       $ 
\begin{array}{lll}
 \hspace {35mm}\widehat{W}(k) & \rightarrow & \ \ \Z \\
 \hspace {35mm} \quad  \downarrow &  & \ \  \downarrow \\
 \hspace {35mm}W(k)& \rightarrow & \Z/2\Z \,
\end{array}
$

\smallskip\n  cf. [Se 03], $\S$ 27.1,  où  l'anneau  $\widehat{W}(k)$ est not\' e $WGr(k)$;
dans ce carr\' e, la fl\`eche horizontale du haut est donn\' ee par le rang, et celle du bas par le rang mod~2. Le fait que le carr\' e soit cart\' esien signifie que l'on peut identifier $\widehat{W}(k)$ au sous-anneau de $W(k)\times \Z$ form\' e des couples $(q,n)$ tels que $n \equiv \rg(q)$ mod $2$.

\medskip

\n {\it Exemples}. (a) Si tout \' el\' ement de $k$ est un carr\' e, on a
$\widehat{W}(k)=\Z$ et $W(k)=\Z/2\Z$.

(b) Si $k=\R$, on a $\widehat{W}(k)\simeq\Z[x]/(x^2-1)$   où $x \in \widehat{W}(k)$ correspond \`a la forme quadratique $\<-1\>$ de rang $1$; on a $W(k)=\Z$, et l'homomorphisme $\widehat{W}(k) \to W(k)$ est donn\' e par $x \mapsto -1$.

\medskip

\n  Les foncteurs $W$ et $\widehat{W}$ ont chacun leurs avantages:

\smallskip

  Les $W(k)$-modules du type $\Inv_k(G)$ (cf. $\S$1.4) sont souvent des modules libres (cf. par exemple prop. 3.2.1 et th. 3.4.1), alors qu'il en est rarement de même des $\widehat{W}(k)$-modules $\Inv_k(G,\widehat{W})$, \`a cause de l'identit\' e $\<1,-1\>.(q - {\rm rang}(q))=~0.$
  
  \smallskip
  
   Les puissances ext\' erieures $\la^i : \widehat{W}(k) \to \widehat{W}(k)$ ([Se 03], 27.1) ne sont pas compatibles (sauf si $-1$ est un carr\' e dans $k$) avec le passage au 
   quotient $\widehat{W}(k) \to W(k)$.
   
   De même, les {\it classes de Stiefel-Whitney} sont d\' efinies sur 
   $\widehat{W}(k)$, mais pas sur $W(k)$. Ce sont des applications $w_i : \widehat{W}(k) \to H^i(k)$, cf. [Se 03], $\S$17.1. Rappelons que, si $q =\<\alpha_1,..., \alpha_n\>$ est une forme quadratique, $w_i(q)$ est la $i$-i\`eme fonction sym\' etrique \' el\' ementaire des \' el\' ements $(\alpha_i)$ de $H^1(k)$. D'apr\`es (1.1.5), appliqu\' e avec $d=1$, on a :
   
   \smallskip (1.2.1) \quad $\sum w_i(q) = \prod_j(1+ w_{2^j})$;
   
   \smallskip \n et: 
   
   \smallskip (1.2.2)  \quad $w_mw_{m'}= e_k^{|A \cap A'|} w_{m\bullet m' }$, \ \  avec les notations de (1.1.8).
   
   \medskip
   
   \n {\it Remarque.}
   
    La formule (1.2.1)  est une cons\' equence imm\' ediate de [Mi 70], Rem. 3.4.
   
 \n   On trouve aussi dans [Mi 70], p. 331, une d\' ecomposition du produit $w_mw_{m'}$ :

\smallskip (1.2.3) \quad $w_mw_{m'}= \sum_i (i, m-i, m'-i) e_k^ iw_{m+m'-i}$,

\smallskip \n   où  $(i, m-i, m'-i) = \frac{(m+m'-i)!}{i!(m-i)!(m'-i)!}$ est le coefficient trinomial associ\' e au triplet 
$\{i,m-i,m'-i\}$, et la somme porte sur les entiers $i$ tels que $0\leqslant i \leqslant \inf(m,m')$. 

\smallskip

Cette formule semble diff\' erente de (1.2.2), puisqu'elle contient une  somme de plusieurs termes. En fait, si l'on utilise les propri\' et\' es de congruence mod 2 des coefficients
trinomiaux \footnote{Si $a = \sum_{\alpha \in A} 2^\alpha, b = \sum_{\beta \in B} 2^\beta, c=\sum_{\gamma \in C} 2^\gamma$, le coefficient trinomial $(a,b,c)= \frac{(a+b+c)!}{a!\cdot b!\cdot c!}$ est impair si et seulement si les ensembles $A, B, C$ sont disjoints, cf. [Di 02].}, on voit que tous ces coefficients sont nuls (dans $\F_2$), \`a la seule exception de celui relatif \`a $i = |A \cap A'| =m+m'-m\bullet m'$, qui est \' egal \`a $1$. Les formules (1.2.2) et (1.2.3) sont donc \' equivalentes.
      
   \smallskip
   
   \n {\it Le cas   où  $-1$ est un carr\' e dans $k$.}
   
    Dans ce cas, les $w_{2i}$ se d\' eduisent de $w_2$ : on a
      $w_{2i} = \gamma_i(w_2)$,   où  $\gamma_i$ d\' esigne la $i$-i\`eme {\it puissance divis\' ee}, qui est d\' efinie sur $H^{\rm{pair}}(k)$, cf. [Vi 09], prop.2.8
      (mais n'est pas d\' efinie sur $H(k)$, contrairement \`a ce qui est affirm\' e dans [Ka 20], App.A).
      Si l'on d\' efinit l'exponentielle
      $\exp(x)$ d'un \' el\' ement $x\in H^2(k)$ comme $1+ \sum_{i\geqslant 1} \gamma_i(x)$, on a donc :
      
\smallskip   (1.2.4)   \quad  $w = w_1\!\cdot\! \exp(w_2)$.

        \medskip     
\n {\bf 1.3. $G$-torseurs}.

\smallskip Rappelons que $G$ d\' esigne un groupe fini.

\smallskip

  Un ensemble $X$, muni d'une action (\`a droite) de $G$ est appel\' e un $G$-{\it torseur} si l'action de $G$ sur $X$ est {\it r\' eguli\`ere}, i.e. est isomorphe \`a celle de $G$ sur lui-même par translations \`a droite. On a $|X| = |G|$.

  Un {\it $G$-torseur sur $k$} est un $k$-sch\' ema \' etale fini $T$  muni d'une action (\`a droite) de $G$ telle que l'ensemble fini $T(k_s)$ soit un $G$-torseur. L'alg\`ebre affine
  de $T$ est une {\it $G$-alg\`ebre galoisienne}, cf. [BS 94], 1.3; c'est une $k$-alg\`ebre \' etale finie de dimension $|G|$. Inversement, le spectre d'une $G$-alg\`ebre galoisienne est
  un $G$-torseur sur $k$.

\smallskip
  Soit $\Tors_G(k)$ l'ensemble des classes
  d'isomorphisme de $G$-torseurs sur $k$. On peut l'identifier aux classes
  d'homomorphismes $\varphi : \Gamma_k \to G$, deux homomorphismes $\varphi_1$ et $\varphi_2$ \' etant dans la même classe s'ils sont $G$-conjugu\' es, cf. [BS 94], $\S$1.3.1. Le torseur trivial  (celui qui a un $k$-point) correspond \`a $\varphi=1$. 
 On peut donc identifier $\Tors_G(k)$ au {\it premier groupe de cohomologie non-ab\' elienne} $H^1(k,G)$, cf. [Se 65], chap. I.5 et [KMRT 98], $S\S$29.

      \medskip
       \n {\bf 1.4. $k$-invariants cohomologiques et invariants de Witt de $G$.}
 
 \smallskip
 
  Soit $\Corps_k$ la cat\' egorie form\' ee des extensions de corps $K/k$, et soit $\Ab$ la cat\' egorie des groupes ab\' eliens. Si $r$ est un entier $\geqslant 0$, notons $H^r$ le foncteur
 $\Corps_k \to \ \Ab$ d\' efini par $K \mapsto H^r(K)$; d\' efinissons de façon analogue les foncteurs $H, \widehat{W}, W, \Tors_G$. 
 
 Un $H^r$\!-$k$-invariant de $G$ est un morphisme
 de foncteurs $ a : \Tors_G \to H^r$. Autrement dit, $a$ est la donn\' ee, pour toute\footnote{On peut se borner aux extensions de $k$ contenues dans une extension alg\' ebriquement close de $k$ de degr\' e de transcendance $\aleph_0$; le cas g\' en\' eral s'en d\' eduit par limite inductive, cf. [Se 03], $\S$1, Remark 1.2.} extension $K/k$, d'une application $a_K : \Tors_G(K) \to H^r(K)$ telle que, pour toute extension $L/K$, le diagramme
 
 \medskip
\hspace{30mm} \ $\Tors_G(K) \quad \stackrel{a_K}{\rightarrow} \ \ H^r(K)$

\hspace{30mm} \quad \quad $\downarrow  \ \hspace{22mm} \downarrow$

\hspace{30mm} \ $\Tors_G(L) \quad \stackrel{a_L}\rightarrow \ \ \ H^r(L)$

\medskip

\n soit commutatif. L'ensemble des $H^r$-$k$-invariants de $G$ est not\' e
$\Inv_k^r(G)$; il a une structure naturelle de $\F_2$-espace vectoriel.

  Dans la d\' efinition ci-dessus  (d\' etaill\' ee dans [Se 03], chap.I), on peut remplacer le foncteur $H^r$ par l'un des trois foncteurs $H, W$,
  $\widehat{W}$. On obtient ainsi des groupes not\' es respectivement $\Inv_k(G), \Inv_k(G,W)$, et
  $\Inv_k(G,\widehat{W})$; leurs \' el\' ements s'appellent des invariants cohomologiques, resp. de Witt, resp. de Grothendieck-Witt, du groupe  $G$ sur $k$. On a $\Inv_k(G) = \boldsymbol{\oplus}_{r\geqslant 0} \Inv_k^r(G)$; ce groupe a une structure naturelle de $H(k)$-alg\`ebre gradu\' ee. De même, $\Inv_k(G,W)$ est une  $W(k)$-alg\`ebre, et $\Inv_k(G,\widehat{W})$ est une $\widehat{W}(k)$-alg\`ebre munie d'une $\la$-structure (au sens de [SGA 6], expos\' e V); ces alg\`ebres sont reli\' ees par un diagramme cart\' esien d\' eduit de celui de 1.2:
  
  \smallskip
  $ 
\begin{array}{lll}
\hspace{30mm}\Inv_k(G,\widehat{W}) & \rightarrow & \  \ \Z \\
 \hspace{30mm}\quad \quad \ \downarrow &  & \ \  \downarrow \quad  \\
 \hspace{29mm} \ \Inv_k(G,W)& \rightarrow & \Z/2\Z \ .
\end{array}
$

  \medskip
  
    Une autre façon de formuler ceci consiste \`a utiliser la d\' ecomposition de $\Inv$
    en somme directe de sa composante constante $\Inv^{\cst}$ et de sa composante normalis\' ee $\Inv^{\norm}$, cf. [Se 03], §4.4 et §4.5; on a :
        
        \medskip
     
     \centerline{ $\Inv_k(G,\widehat{W})^{\cst} = \widehat{W}(k)$ \ \ et \ \ $\Inv_k(G,W)^{\cst} = W(k)$,}

\medskip
     \centerline{ $\Inv_k(G,\widehat{W})^{\norm} = \ \Inv_k(G,W)^{\norm}$.}

  \medskip

\n  Dans la suite, nous utiliserons en g\' en\' eral $\Inv_k(G,W)$, et nous laisserons au lecteur le soin d'en d\' eduire les r\' esultats correspondants pour $\Inv_k(G,\widehat{W})$.
  
  \medskip

  \n{\bf 1.5. Fonctorialit\' e en $k$.}
  
  \smallskip   Soit $\cc$ l'un des foncteurs $H, W, \widehat{W}$. Si $k'$ est une extension de $k$, tout \' el\' ement
    de $\Inv_k(G,\cc)$ d\' efinit par restriction un invariant de $\Inv_{k'}(G,\cc)$; d'où un homomorphisme $\Inv_k(G,\cc) \to \Inv_{k'}(G,\cc)$.
    
    \smallskip
    
     \n {\sc Proposition 1.5.1}. {\it Soit $k'$ une extension finie de $k$ de degr\' e impair. L'application $\Inv_k(G,\cc) \to \Inv_{k'}(G,\cc)$ est injective.}
      
      \medskip  
      \n {\it D\' emonstration}. 
  
  \smallskip 
  
  \n {\sc Lemme 1.5.2}. {\it Soit $K$ une extension de $k$. Il existe un corps quotient $K'$ de $K \otimes_k k'$ tel que $[K':K]$ soit impair}.
  
  [Autrement dit, il existe une extension compos\' ee de $K$ et $k'$ ([AV], $\S$2.4) dont le degr\' e sur $K$ est impair.]
  
  \smallskip

  \n {\it D\' emonstration du lemme 1.5.2}. Soit $k''$ la plus grande extension s\' eparable de $k$ contenue
  dans $k'$. Comme la $k$-alg\`ebre $k''$ est \' etale de degr\' e impair, il en est de même de la $K$-alg\`ebre $A= K \otimes_k k''$; cette $K$-alg\`ebre est donc un produit de corps; l'un d'eux, disons $K''$, est de degr\' e impair sur $K$.
  
  L'alg\`ebre $A'= K''\otimes_{k''}k'$ est un quotient de $A\otimes_{k''}k' = K \otimes_kk'$. Soit  $K'$ un corps quotient de $A'$; c'est un quotient de $K \otimes_k k'$. Comme  $k'$ est une extension radicielle de $k''$ ([A IV-VII], $\S$7, prop. 13), $K'$ est une extension radicielle de $K''$,
  donc de degr\' e impair sur $K''$, et donc aussi de degr\' e impair sur $K$.
  
  \smallskip
  
  \n {\it Fin de la d\' emonstration de la prop. 1.5.1}. Soit $a \in \Inv_k(G,\cc)$, et soit $a'$
  son image dans $\Inv_{k'}(G,\cc)$. Supposons que $a'=0$ et montrons que $ a=0$. Soit $T$
  un $G$-torseur sur une extension $K$ de $k$. On a $a(T) \in \cc(K)$. Soit $K'$
  une extension de $K$ ayant la propri\' et\' e du lemme 1.5.2, et soit $T'$ le $G$-torseur
  sur $K'$ d\' eduit de $T$ par extension des scalaires. Puisque $K'$ est une extension de $k'$, on a $a(T') = a'(T')=0.$ D'autre part, $a(T') \in \cc(K')$ est l'image de $a(T) \in \cc(K)$ par l'application $\cc(K) \to \cc(K')$, application qui est injective puisque $[K':K]$ est impair : lorsque $\cc = W$ ou $ \widehat{W}$, c'est un r\' esultat de Springer ([Sp 52]) - le cas $\cc = H$ est bien connu. On a donc $a(T)=0$, ce qui ach\`eve la d\' emonstration.
  
  \medskip

   \n {\bf 1.6. Fonctorialit\' e en $G :$ induction et restriction.}
   
   \smallskip
 
 Soit $f:G' \to G$ un homomorphisme de groupes finis. Il lui correspond une application d'{\it induction}
 $f_*:\Tors_{G'}(K) \to \Tors_G(K)$ pour tout $K/k$; cette application associe \`a un
 $G'$-torseur $T'$ sur $K$ le  $G$-torseur $f_*T'$ sur $K$ qui s'en d\' eduit au moyen de $f$, \`a savoir le quotient de $T'\times G$ par $G'$ agissant par	$g'.(t',g) =
 (t'g',f(g')^{-1}g)$. 
 Avec les notations de [Se 65], $\S$I.5.3, on a  $$f_*T' = T' \times^{G'}G.$$ 
 
  Si $T'$ est d\' efini par un homomorphisme $\varphi:\Gamma_K \to G'$,
 alors $f_*T'$ est d\' efini par le compos\' e $f \circ \phi : \Gamma_K \to G' \to G$.

 \smallskip \n{\sc Proposition 1.6.1}. {\it Supposons que $f : G' \to G$ soit injectif.
 Soit $T$ un $G$-torseur sur $K$. Pour qu'il existe un $G'$-torseur $T'$ tel que
 $T' \simeq f_*T'$, il faut et il suffit que le $K$-sch\' ema $T/G'$, quotient de $T$ par l'action de $G'$, ait un $K$-point.}
 
   \smallskip Cela r\' esulte de [Se 65], $\S$I.5.4, prop. 37.
   
   \smallskip
 
 Soit $\cc$ comme dans 1.5. Si $a \in \Inv_k(G,\cc)$, on d\' efinit $f^*(a) \in \Inv_k(G',\cc)$ par la formule $f^*(a)(T')= a(f_*T')$. D'où un homomorphisme
 
 \smallskip
 \centerline{$f^* : \Inv_k(G,\cc) \to \Inv_k(G',\cc).$}
 
 \smallskip
   Lorsque $G'$ est un sous-groupe de $G$, et que $f$ est l'injection de $G'$ dans $G$,
   l'homomorphisme $f^*$ s'appelle la {\it restriction} de $G$ \`a $G'$; on le note
   $\Res^{G'}_G$, ou simplement $\Res$.
      
      \smallskip
   
   \n {\sc Proposition 1.6.2}. {\it Si $G'$ est un sous-groupe de $G$ d'indice impair, l'homomorphisme $\Res: \Inv_k(G,\cc) \to \Inv_k(G',\cc)$ est injectif.}
   
   \smallskip
   
   \n {\it D\' emonstration}. Le cas  $\cc = H$ est trait\' e dans [Se 03], $\S$14 grâce \`a  la construction d'un homomorphisme de {\it corestriction} $\Cor : \Inv_k(G',H) \to \Inv_k(G,H)$ tel que
   $\Cor \circ \Res (x)= (G:G')\!\cdot \!x$ pour tout $x \in \Inv_k(G)$. On peut appliquer la même m\' ethode aux foncteurs $W$ et $\widehat{W}$,
   mais cela demande quelques modifications: la plus importante est que $\Cor \circ \Res$, au lieu d'être la multiplication par $(G:G')$, est la 
   multiplication par l'\' el\' ement $\Cor(1)$ de $\widehat{W}(k)$, qui est de rang $(G:G')$, donc impair, ce qui entraîne que $\Cor \circ \Res$ est
   injectif, cf. [La 05], VIII.8.5). 
  
  \medskip
   \n Pour la commodit\' e du lecteur, voici une d\' emonstration directe.
      
    Soit $a \in \Inv_k(G,\cc)$, et soit $a' = \Res(a) \in \Inv_k(G',\cc)$. Supposons que $a'=0$, et montrons que $a=0$. Si $T$ est un $G$-torseur sur une extension $K/k$, on doit prouver que l'\' el\' ement $a(T)$ de $\cc(K)$ est nul. C'est clair si $T$ est l'induit d'un $G'$-torseur $T'$ car on a alors $a(T)=a'(T')$. On va se ramener \`a ce cas:
     
      Le $K$-sch\' ema $T/G'$ est fini \' etale, donc 
  de la  forme $\prod \Spec(K_i)$,  où les $K_i$ sont des extensions finies s\' eparables de $K$; on a
  $\sum [K_i:K] = (G:G')$, qui est impair. L'un des $K_i$ est donc de degr\' e impair. 
  Choisissons un tel $K_i$. Comme $T/G'$ a un $K_i$-point, la prop. 1.6.1 montre que $T$ devient l'induite d'un $G'$-torseur
  apr\`es extension du corps \`a $K_i$. L'image de $a(T)$ dans $\cc(K_i)$ est donc  $0$.  Comme $[K_i:K]$ est impair, cela entraîne $a(T)$=0.

  \medskip
  
  \n{\bf 1.7. Sous-groupe normal.}
  
  \smallskip  
  
    Soit $1 \to A \to B \to C \to 1$ une suite exacte de groupes finis. Le groupe $B$
    agit sur $A$ par $b(a)= bab^{-1}$. Cela donne une action de $B$ sur $\Inv_K(A,\cc)$, donc aussi une action de $A$. 
    
    \smallskip\n {\bf Proposition 1.7.1.} {\it L'action de $A$ sur $\Inv_K(A,\cc)$ est triviale.}
  
    \smallskip Cela r\' esulte de ce que l'action de $A$ sur $\Tors_K(A)$ est triviale.
    
    \smallskip On a donc une action naturelle du groupe $C=B/A$ sur $\Inv_K(A,\cc)$.
    Cette action fixe l'image de $\Res : \Inv_K(B,\cc) \to \Inv_K(A,\cc)$.
    
     \medskip     
  \n {\bf 1.8. Bases g\' en\' eriques et invariants d'un produit de deux groupes.}

  \smallskip

    Soit $\cc$ l'un des foncteurs $H$ et $W$. Soit $(a_j)_{j\in J}$ une famille finie d'\' el\' ements de $\Inv_k(G,\cc)$. Nous dirons que cette famille est une {\it base g\' en\' erique} si elle a la propri\' et\' e suivante :
    
    \smallskip
    
   {\it Pour toute extension $K/k$, les images des $a_j$ dans $\Inv_K(G,\cc)$
     forment une base du $\cc(K)$-module $\Inv_k(G,\cc)$.}
     
     \smallskip
   \n  [Dans [Ga 24], une telle base est appelée ``strong basis''.]
   
   \smallskip
     
     \n {\it Exemple.} 
     
     Si $G = S_n$ et $\cc = H$ (resp. $W$), les classes de Stiefel-Whitney $w_i(q)$ (resp. les
     $\la^iq$) avec $0 \leqslant i \leqslant [\frac{n}{2}]$) forment une base g\' en\' erique de $\Inv_k(G,\cc)$, cf. $\S$3.1.
     
     \smallskip
     
     L'existence d'une base g\' en\' erique pour un groupe $G$ permet de ramener
   les invariants d'un produit $G \times G'$ \`a ceux de $G$ et de $G'$. De façon plus pr\' ecise:
     
     \smallskip
     \n {\sc Proposition 1.8.1.} {\it Soit $(a_j)_{j\in J}$ une base g\' en\' erique de $\Inv_k(G,\cc)$. Soit $G'$ un groupe fini, et soit $a \in \Inv_k(G\times G',\cc)$. Il existe alors une famille unique d'\' el\' ements $b_j \in 
     \Inv_k(G',\cc)$ telle que l'on ait $$ a(T\times T')= \sum_{j\in J} a_j(T)b_j(T')$$ quels que
     soient l'extension $K/k$, le $G$-torseur $T$ sur $K$ et le $G'$-torseur $T'$ sur~$K$.}
     
     \n [Noter que $\Tors_K(G \times G') = \Tors_K(G)\times \Tors_K(G')$.]
     
     \smallskip
     \n {\it D\' emonstration} (cf. [Se 03, exerc.16.15). Soit $T'$ un $G'$-torseur sur une extension $K/k$. Pour tout $G$-torseur $T$ sur une extension $L$ de $K$,
     soit $b(T) = a(T\times T') \in \cc(L)$. L'application $T \mapsto b(T)$ est un
       $\cc$-invariant de $G$ sur $K$. D'apr\`es (BG), on peut l'\' ecrire de façon unique 
       sous la forme $b(T)= \sum_{j\in J} a_j(T)b_j(T')$ avec $b_j(T') \in \cc(K)$. Pour tout $j$, l'application $T' \mapsto b_j(T')$ est un $\cc$-invariant de $G'$ sur $k$.
       D'où la proposition.
       
       \smallskip
       
 \n {\it Remarque.}  Sans faire d'hypoth\`ese sur $G$, ni sur $G'$, on a une
       application naturelle
       
       \smallskip
       
   \quad $\Inv_k(G,\cc) \otimes_{\cc(k)} \Inv_k(G',\cc) \ \to \Inv_k(G\times G',\cc),$
      
      \smallskip
      
 \n qui transforme un produit tensoriel $a \otimes a'$, $(a\in \Inv_k(G,\cc), a'\in \Inv_k(G',\cc)$) en le $G\times G'$-invariant $(T,T') \mapsto a(T)\!\cdot\! a'(T')$.
 
    La proposition 1.8.1
   \' equivaut \`a dire que {\it cette application est un isomorphisme  si $G$ a une base g\' en\' erique}.
   En particulier:
   
   \smallskip
   
   \n {\sc Corollaire 1.8.2.} {\it Si  $G'$ a une base g\' en\' erique $(a'_{j'})_{j'\in J'}$,
   alors $G\times G'$ a pour base g\' en\' erique $(a_j\!\cdot\! a'_{j'})_{(j,j')\in J\times J'}$.}
   
   \smallskip
   Par r\' ecurrence, on en d\' eduit :   
   
   \smallskip
   
   \n {\sc Proposition 1.8.3}. {\it Soient $(G_1,$\dots$,G_m)$ des groupes
   finis ayant chacun une base g\' en\' erique pour $\cc$. Il en est alors de même de $G = \prod G_i$,  et l'on a $:$
   
   \smallskip
   
  \centerline{$\Inv_k(G) \ \ = \ \ \Inv_k(G_1,\cc)\  \otimes \ \Inv_k(G_2,\cc) \ \otimes \ \cdots \ \otimes \ \Inv_k(G_m,\cc),$}
  
  \smallskip
   
  \n  où  les produits tensoriels sont relatifs \`a l'anneau} $\cc(k)$.
            
  \bigskip
  \n {\sc $\S$2. Relations entre invariants de Witt et invariants cohomologiques}
    
      \medskip
    
       \n {\bf 2.1. La filtration de l'anneau de Witt d'un corps.}
       
       \smallskip
       
       Soit $K$ un corps de caract\' eristique $\ne 2$, et soit $I_K$ le noyau de  l'homomorphisme $\rg : W(K) \to \Z/2\Z$. Les puissances 
       $I_K^n$ de $I_K$ d\' efinissent une filtration descendante sur $W(K)$. 
        Si $x \in W(K)$, on note $\fil(x)$ la borne inf\' erieure des $n$ tels que $x\in I^n_K$; d'apr\`es un th\' eor\`eme d'Arason-Pfister ([AP 71]), on a $ \fil(x)=\infty$ si et seulement si $x=0$.
        
        \smallskip
       Soit $n$ un entier $\geqslant0$; notons $\gr^n W(K)$ le quotient $I^n_K/I^{n+1}_K$;
       c'est un $\F_2$-espace vectoriel.
D'apr\`es les conjectures de Milnor ([Mi 70]), rappel\' ees au $\S$1.1, le groupe $\gr^n W(K)$ est canoniquement isomorphe \`a $H^n(K)$. On obtient ainsi un isomorphisme entre les $\F_2$-alg\`ebres gradu\' ees
$H(K)$ et       $\gr W(K) = \sum_n\gr^n W(K)$. Cet isomorphisme associe \' eà
un \' el\' ement d\' ecomposable $(x_1)\cdots(x_n)$ de $H^n(K)$ la classe de la $n$-forme de Pfister $\<1,-x_1\> \<1,-x_2\>\cdots\<1,-x_n\>$.

\medskip

\n{\small[Plus pr\' ecis\' ement, notons $\Mil(K)$ l'alg\`ebre de Milnor de $K$. Dans [Mi 70], Milnor d\' efinit des homomorphismes:

\   $\alpha: \Mil(K)/2\Mil(K) \to H(K)$ \ \ et \ \  $\beta : \Mil(K)/2\Mil(K) \to \gr W(K)$,

\n qui sont des isomorphismes d'apr\`es [Vo 03.I], [Vo 03.II], [OVV 07], [Mo 05]. L'isomorphisme utilis\' e ici est $\beta \circ \alpha^{-1}$.]}

    \medskip
    
    \n {\bf 2.2. La filtration de $\Inv_k(G,W)$.}
    
    \smallskip
     Soit  $a \in \Inv_k(G,W)$. Soit $\fil(a)$ la borne inf\' erieure (finie ou infinie)
    des $\fil(a(T))$, pour tous les $G$-torseurs $T$ sur toutes les extensions $K/k$. On a
    $\fil(a) = \infty$ si et seulement si $a=0$. Si $n \geqslant 0$, notons $\Inv_k(G,W)_n$ le sous-groupe de $\Inv_k(G,W)$ form\' e des \' el\' ements $a$ tels que $\fil(a) \geqslant n$. On obtient ainsi une filtration sur l'anneau $\Inv_k(G,W)$; soit $\gr \Inv_k(G,W)$ le gradu\' e associ\' e; on a:
    
 \smallskip  \centerline{  $\gr^n \Inv_k(G,W) = \Inv_k(G,W)_n/\Inv_k(G,W)_{n+1}$.}
        
     \smallskip   Soit $a \in \Inv_k(G,W)_n$, et soit $T$ un
 $G$-torseur sur une extension $K$ de $k$. On a $a(T) \in I^n_K$. 
    Soit $\tilde{a}(T)$ l'image de $a(T)$ dans $H^n(K) \simeq I^n_K/I^{n+1}_K$. On obtient ainsi un invariant cohomologique $\tilde{a} \in \Inv_k^n(G)$; cet invariant est nul si et seulement si $\fil(a) > n$. D'où un homomorphisme :
    
    \smallskip
    
    \centerline{$\gr^n\Inv_k(G,W) \to \Inv_k^n(G).$}

    \medskip
    
    \n {\sc Proposition 2.2.1.} {\it L'homomorphisme  $\gr^n\Inv_k(G,W) \to \Inv_k^n(G)$ d\' efini ci-dessus est injectif.}   
    
    C'est clair.
    
    \smallskip
    
  \n  On obtient ainsi un plongement de $\gr\Inv_k(G,W)$ dans $\Inv_k(G)$. 
  
  \smallskip     
     {\it Remarque}. Les résultats des §§ 2.1 et 2.2 ne sont pas nouveaux; on les trouve dans [Ga 20], et, sous une forme un peu différente, dans [Hi 09] et [Hi 20].     
              
    \medskip
    
    \n {\bf 2.3. Application aux homomorphismes de restriction.}
    
    \smallskip
    Soit $G'$ un sous-groupe de $G$.
    
    \smallskip
    
    \n {\sc Proposition 2.3.1.}{\it \ Supposons que \ $\Res_H : \Inv_k(G) \to \Inv_k(G')$ soit injectif $;$  alors $\Res_W : \Inv_k(G,W) \to \Inv_k(G',W)$ est injectif.}
    
    \smallskip
    \n [Autrement dit, tout th\' eor\`eme de d\' etection pour les
    invariants cohomologiques mod  2  entraîne un th\' eor\`eme analogue pour les invariants de Witt.]

    \smallskip
   \n {\sc Lemme 2.3.2.} {\it Soit $a \in \Inv_k(G,W)$ et soit $a' = \Res_W(a) \in \Inv_k(G',W)$.
   Alors $\fil(a') = \fil(a).$}
   
   \smallskip
   
   \n {\it D\' emonstration du lemme}. Soit $d = \fil(a)$. On peut supposer $a\ne 0$, i.e. $d \ne \infty$.    Il est imm\' ediat que $\fil(a') \geqslant d$. Soit $\alpha$ l'image de $a$
   dans  $\gr^n\Inv_k(G,W) \simeq \Inv_k^n(G)$ et soit $\alpha'$ l'image de $a'$ dans
   $\Inv_k^n(G')$; on a  $\alpha' = \Res_H(\alpha)$. Comme $\fil(a)= d$, on a $\alpha \ne 0$; puisque $\Res_H$ est injectif, on a $\alpha' \ne 0$, ce qui \' equivaut \`a $\fil(a')=d$.
   
   \smallskip
   
  \n {\it D\' emonstration de la prop. 2.3.1}.  Le lemme 2.3.2, appliqu\' e \`a un \' el\' ement $a$ du noyau de $\Res_W(G,W)$, montre que $\fil(a) = \fil(0) = \infty$, d'où $a=0$.
    
    \medskip
    \n  {\it Remarque}.
Les r\' esultats et les d\' emonstrations de ce $\S$ s'\' etendent au cas   où  le groupe $G$, au lieu d'être un groupe fini, est un $k$-groupe alg\' ebrique lin\' eaire.
         
       \bigskip

     \n {\sc $\S$3. Invariants de $A_n$: pr\' eliminaires}

  \medskip
  
   \n {\bf 3.1. Rappels sur les torseurs et les invariants de $ S_n$.}
  
  \smallskip
  
    Dans ce qui suit, $n$ d\' esigne un entier $>1$.    Soit $S_n$ le groupe des permutations de $\{1,$\dots$,n\}$. L'ensemble $\Tors_{S_n}\!(k)$
  s'identifie, via la th\' eorie de Galois, \`a l'ensemble $\Et_n(k)$ des classes d'isomorphisme de $k$-alg\`ebres \' etales de dimension $n$, cf. [A IV-VII], p. A V 72, et [Se 03], §3.2. Si $L$ est une telle alg\`ebre, on note $q_L$ sa forme trace, autrement dit la forme quadratique
de rang $n$ qui  associe \`a un \' el\' ement $x$ de $L$ l'\' el\' ement $\Tr_{L/k}(x^2)$ de $k$. Cela d\' efinit une application  $\Tors_{S_n}(k) \to \widehat{W}(k)$. En appliquant cette construction
aux extensions $K$ de $k$, on obtient un invariant de Grothendieck-Witt de $S_n$,
not\' e $q$. Ses puissances ext\' erieures $\la^ i(q)$ sont aussi
des \' el\' ements de $\Inv_k(S_n,\widehat{W})$; soient $\la^iq$ leurs images dans
$\Inv_k(S_n,W)$, et soient  $w_i(q) \in  \Inv^i_k(S_n)$ les classes de Stiefel-Whitney 
de $q$. D'apr\`es [Se 03], §25.13 et §29.2, on a :

\smallskip
   
\n {\sc Proposition 3.2.1}. {\it Les $w_i(q)$ avec   $0 \leqslant i \leqslant  [\frac{n}{2}] $  forment une base du $H(k)$-module} $ \Inv_k(S_n)$.

\medskip

\n {\sc Proposition 3.2.2}. {\it Les $\la^iq$ avec $0 \leqslant i \leqslant  [\frac{n}{2}] $ forment une base du 
$W(k)$-module $\Inv_k(S_n,W)$}.

    \medskip
   
   \n {\it Remarques.} 
   
   1. Lorsqu'il est utile de mentionner explicitement $n$, ou le couple $(n,k)$, on \' ecrit $q_n$ ou $q_{n,k}$ \`a la place de $q$. Par exemple, lorsqu'on compare $q$ \`a sa restriction
   au sous-groupe $S_{n-1}$ de $S_n$, on \' ecrit $ \Res(q_n) = q_{n-1}+1 $. 
     
    2. Le base donn\' ee par la prop. 3.2.1 est une base g\' en\' erique, au sens du  §1.8: cela  se voit en appliquant la proposition aux extensions de $k$. Il en est de même de la base de la prop. 3.2.2, ainsi que de celles du th. 3.4.1 et du th. 3.4.2.
    
    3. Dans [Se 03], §25.3, il y a des classes de Stiefel-Whitney ``galoisiennes'', not\' ees  $w_i^{\gal}$, qui sont li\' ees aux $w_i(q)$ par la formule suivante, due \`a B. Kahn ([Ka 84], voir aussi [Se 03], version corrigée, §25.10):

   \smallskip
   
$   (3.2.1)      \  w_i^{{\rm gal}} = \left\{
\begin{array}{lll}
 w_i(q) + (2)_k\!\cdot\! w_{i-1}(q) \ \ \ {\rm si} \ i  \ {\rm est \ pair}\\
w_i (q)\  \hspace{23mm} {\rm si} \ i \ {\rm est \ impair},\ \\
\end{array}
\right.
$ 

\n \smallskip ou, ce qui revient au même :

\smallskip (3.2.2)  $w^{\gal} = w\!\cdot\!(1+ (2)_k\!\cdot\! w_1)$, \   où  \ $w^{\gal} = \sum_i w_i^{\gal}.$

  \smallskip \n  On a :
 
 \smallskip  (3.2.3) $w^{\gal} = \prod_j(1 + w^{\gal}_{2^j})$.
 
  \smallskip 
  Cela se d\' emontre en utilisant  (1.2.1) et (3.2.2), et c'est un cas particulier d'un r\' esultat plus g\' en\' eral de [Ka 84].
  
  \smallskip
  
La formule (3.2.1) montre que les $w_i^{{\rm gal}}$,
    avec $0 \leqslant i \leqslant [\frac{n}{2}]$, forment, elles aussi, une base de $\Inv_k(S_n)$.   
    
    \medskip
    
    \n {\bf 3.3. $A_n$-torseurs.}
    
    \medskip
    
     Soit $A_n$ le sous-groupe de $S_n$ form\' e des \' el\' ements pairs (i.e. de signature 1).  La suite exacte $1 \to A_n \to S_n \to \{\pm 1\} \to 1$ (valable car $n >1$) montre que tout torseur 
      sous $A_n$ d\' etermine un torseur sous $S_n$ dont l'image par 
      
      \smallskip \centerline {$ \Tors_k(S_n ) \ \to  \ \Tors_k( \{\pm 1\})\ \simeq \ H^1(k,\Z/2\Z)$}
      
      \smallskip \n  est 0. En d'autres termes, un $A_n$-torseur $t$ d\' efinit une
      $k$-alg\`ebre \' etale $L_t$ telle que $q_{L_t}$ soit de discriminant 1
      (dans $k^\times/k^{\times  2}$), i.e. $w_1(q_{L_t})=0$.
      
      R\' eciproquement,  une $k$-alg\`ebre \' etale $L$ de dimension $n$ telle que $w_1(q_L)=0$ provient de {\it un ou deux}
      $A_n$-torseurs. En effet, $L$ est donn\' ee par une classe de $S_n$-conjugaison d'homomorphismes de $ \Gamma_k $ dans $ S_n$, et l'hypoth\`ese $w_1(q_L)=0$ signifie
      que les images de ces homomorphismes sont contenues dans $A_n$. 
     Si ces homomorphismes sont  $A_n$-conjugu\' es, ils d\' efinissent un $A_n$-torseur unique $t$ tel que $L \simeq L_t$. 
     Sinon,  ils forment deux classes de $A_n$-conjugaison qui donnent deux 
      $A_n$-torseurs $t, t'$ tels que $L \simeq L_t \simeq L_{t'}$; ces deux torseurs ne sont pas isomorphes, mais ils sont transform\' es l'un de l'autre par l'action de $S_n/A_n$ sur $\Tors_k(A_n)$, cf. 1.7.
      
     {\small Une façon de distinguer $t$ et $t'$ consiste \`a {\it orienter} $L$, c'est-\' eà-dire \`a  choisir un \' el\' ement  $e$ de $ \wedge^nL$ tel que $e\!\cdot \!e = 1$, le produit scalaire \' etant celui d\' efini par $\wedge^nq_L$. Il y a deux choix possibles, qui correspondent \`a $t$ et $t'$. Lorsque $L$ est donn\' ee par une \' equation polynomiale $F(x)=0$,   où  $F$ est de degr\' e $n$, l'hypoth\`ese $w_1(q)=0$ \' equivaut \`a dire que le discriminant de $F$ est un carr\' e, et le choix de  $t$ ou $t'$ se traduit par le choix d'une racine carr\' ee de ce discriminant.}
    
    \medskip

   \n {\bf 3.4. Invariants de $ A_n$} : {\bf \' enonc\' es des r\' esultats.}
  
   \smallskip

\n L'injection de $A_n$ dans $S_n$ d\' efinit par restriction des homomorphismes

\smallskip

 \quad  $ \Res :\Inv_k(S_n,\cc) \ \ \to \ \ \Inv_k(A_n,\cc)$ \ \   où  $\cc$ est l'un des foncteurs $H, W,\widehat{W}$.

\smallskip
  
\n  L'invariant $q \in \Inv_k(S_n,\widehat{W})$ d\' efinit un \' el\' ement
  de $\Inv_k(A_n,\widehat{W})$ que nous noterons  $q^a$ (ou $q^a_n
  $, si n\' ecessaire). Comme dans 3.1, on en d\' eduit des invariants 
  
  \smallskip 
  
\centerline{  $w_i(q^a) \in \Inv_k^i(A_n), \ \ \la^ i(q^a) \in \Inv_k(A_n,\widehat{W})$ \ et \ $ \la^iq^a \in \Inv_k(A_n, W)$.} 
    
    \smallskip
    
    Ces invariants sont nuls pour $i>n$. On a $w_1(q^a)=0$ et $\la^n(q^a)=1$. Cette derni\`ere formule entraîne que $ \la^i(q^a) = \la^{n-i}(q^a)$ pour tout $i$.
  \smallskip
  
  \n Les analogues pour $A_n$ des prop. 3.2.1 et prop. 3.2.2  sont : 
  \medskip
     
\n {\sc Th\' eor\`eme 3.4.1}. {\it Les $w_{2i}(q^a)$ avec $0 \leqslant i \leqslant  [\frac{n}{4}] $ forment une base du $H(k)$-module $\Inv_k(A_n)$}.

\smallskip On a $w_i(q^a)= 0$ si $i$ est impair ou si $i > 2[\frac{n}{4}]$: le premier cas r\' esulte 
de la formule $w_{2j+1}= w_1w_{2j}$ et le second cas r\' esulte du premier et du th. 25.13 de [Se 03].

\medskip
  
\n {\sc Th\' eor\`eme 3.4.2}. {\it Les $\la^iq^a$ avec $0 \leqslant i \leqslant  [\frac{n}{4}] $ forment une base du 
$W(k)$-module $\Inv_k(A_n,W)$}.

\smallskip

En particulier, les $\la^iq^a, \ i > [\frac{n}{4}]$, sont des combinaisons $W(k)$-lin\' eaires des
$\la^iq^a$, avec $i \leqslant \frac{n}{4}$; nous verrons plus loin (th. 6.1.1) qu'ils en sont même des {\it combinaisons $\Z$-lin\' eaires \`a coefficients ind\' ependants de } $k$.

\smallskip

\n {\it Remarques.}

1. On a $w_i(q^a)= 0$ si $i$ est impair ou si $i > 2[\frac{n}{4}]$: le premier cas r\' esulte 
de la formule $w_{2j+1}= w_1w_{2j}$ et le second cas r\' esulte du premier et du th. 25.13 de [Se 03]. 

2. Les classes galoisiennes $w_i^{{\rm gal}}$ mentionn\' ees \`a la fin du $\S$3.1 sont \' egales aux classes $w_ i$ : cela r\' esulte de (3.2.1).
 
 \smallskip
 
 \n {\sc Corollaire 3.4.3}. {\it L' homomorphisme $\Res :  \Inv_k(S_n,\cc) \ \ \to \ \ \Inv_k(A_n,\cc)$  est surjectif.}
 
 \n [Rappelons que $\cc$ d\' esigne l'un des foncteurs $H, W,\widehat{W}$.]
 
 \smallskip  \n {\it D\' emonstration.} Les $w_i(q^a)$ et les $\la^i(q^a)$ sont les restrictions des $S_n$-invariants $w_i(q)$ et $\la^i(q)$; cela entraîne le corollaire lorsque $\cc = H$ ou $W$; le cas de $\widehat{W}$ r\' esulte de celui de $W$.

 \smallskip \n {\sc Corollaire 3.4.4}. {\it Soient $t,t'$ deux $A_n$-torseurs sur $k$ qui d\' efinissent des $k$-alg\`ebres \' etales isomorphes}  (cf.  $\S$3.3). {\it Alors $a(t)=a(t')$ pour tout
 $a \in \Inv_k(A_n,\cc)$.}
 
 \smallskip Cela r\' esulte du cor.  3.4.3.
    
     \medskip
     
     Les d\' emonstrations du th. 3.4.1 et du th. 3.4.2 seront donn\' ees au §5.5. La m\' ethode est semblable \`a celle
   utilis\' ee pour $S_n$ dans [Se 03]; elle repose sur la construction
   d'un sous-groupe $E$ de $A_n$ ayant les deux propri\' et\' es suivantes :
     
     \smallskip
  $  \bullet$   Ses $\cc$-invariants sont connus (ou faciles \`a d\' eterminer);
      
      \smallskip
   $ \bullet$  L'homomorphisme $\Res: \Inv_k(A_n,\cc) \to \Inv_k(E,\cc)$ est injectif,
      autrement dit $E$ ``d\' etecte'' les $\cc$-invariants de $A_n$.  
      
      \smallskip     
     
  Comme on le verra au $\S$4, on peut prendre pour $E$ un sous-groupe ab\' elien de $A_n$ engendr\' e par des bitranspositions. 
  
\n  [Nous appelons {\it bitransposition} un produit de deux transpositions \`a supports disjoints, autrement dit une involution qui fixe tous les points sauf quatre. L'ensemble des points qui ne sont pas fix\' es est appel\' e le {\it support} de la bitransposition.  Toute partie \`a quatre \' el\' ements est le support de trois bitranspositions qui commutent entre elles.]     
         
    \bigskip
   
    \n {\sc $\S$4. Sous-groupes de $A_n$ d\' etectant les invariants}
    
    \medskip
  
  \n  {\bf 4.1. D\' efinition des sous-groupes $E$ et $A$.}
  
  \smallskip 
    Soit $n$ un entier $>0$, soit $m=[\frac{n}{4}]$, et soit $c=n-4m$. Soit $X  =\{1,$\dots$,n\}$. Pour $i= 1,\dots,m$, soit $X_i=\{4i-3, 4i-2,4i-1,4i\}$, et soit $X_0$ l'ensemble $\{4m+1,\dots,4m+c\}$. La famille  $(X_0,X_1,\dots, X_m)$ est une partition de $X$.
    
Soit $A(i)$ le sous-groupe de $A_n$ form\' e des \' el\' ements qui fixent $X\sm X_i$;
c'est un groupe isomorphe \`a $A_4$. Soit $E(i)$ le 2- sous-groupe de Sylow de $A(i)$; c'est un groupe ab\' elien \' el\' ementaire d'ordre 4, dont les \' el\' ements $\ne 1$ sont les bitranspositions
de support $X_i$. Soient $E$  et $A$ les sous-groupes de $A_n$ d\' efinis par :
    $$E = \prod E(i) \quad {\rm et} \quad A = \prod A(i).$$
    Le groupe $A$ est isomorphe \`a $A_4 \times \cdots \times A_4$ ($m$ facteurs).
    Le groupe $E$ est l'unique $2$-Sylow de $A$; c'est un 2-groupe ab\' elien \' el\' ementaire de rang $2m$. Les groupes $A$ et $E$ fixent
 $X_0$.
    
    \smallskip
    \n {\it Exemple.} Si $n=10$, on a $m=2$, $\rg(E)=4$ et l'on peut prendre pour base de $E$ les quatre bitranspositions :  $(1 \ 2)(3 \ 4),  \ (1 \ 3)(2 \ 4), \  (5 \ 6)(7 \ 8),  \ (5 \ 7)(6 \ 8)$. On a $X_0= \{9,10\}.$
    
    \smallskip \n {\it Remarque}. A conjugaison pr\`es, $E$ est caract\' eris\' e par les propri\' et\' es suivantes :
    
     (a) C'est un $2$-sous-groupe ab\' elien \' el\' ementaire de rang $2m$, engendr\' e
     par des bitranspositions.
     
     (b) Si  une bitransposition $s$ appartient  \`a $E$, alors les deux bitranspositions de même support que $s$ appartiennent aussi \`a $E$.
     
         \medskip
    
    \n {\bf 4.2. Le th\' eor\`eme de d\' etection.}
    
    \smallskip
    
    Ce th\' eor\`eme dit que le sous-groupe $E$ d\' efini au $\S$4.1 d\' etecte les $\cc$-invariants de $A_n$, si $\cc = H, W$ ou $\widehat{W}$:
    
    \smallskip

 \n {\sc Th\' eor\`eme 4.2.1}. {\it L'homomorphisme de restriction $\Inv_k(A_n,\cc) \to
 \Inv_k(E,\cc)$ est injectif.}
 
 \smallskip Comme $A$ contient $E$, on a :
 
 \smallskip
 
 \n {\sc Corollaire 4.2.2.} {\it  L'homomorphisme de restriction $\Inv_k(A_n,\cc) \to
 \Inv_k(A,\cc)$ est injectif.}
 
         \smallskip        \n {\sc Corollaire 4.2.3}. {\it Pour tout $m\geqslant 0$, les applications de restriction $:$
   
   \smallskip      
        $\Inv_k(A_{4m+3},\cc) \to \Inv_k(A_{4m+2},\cc) \to  \Inv_k(A_{4m+1}, \cc) \to \Inv_k(A_{4m},\cc)$
       
        \smallskip
          \n      sont injectives.} 
         \n {\small [On verra au §5.5 que ces applications sont même bijectives.]}
       
       \smallskip
        Cela r\' esulte du fait que le sous-groupe $E$ est le même pour
        $A_{4m}, $\dots$, A_{4m+3}$.

    \medskip
        \n {\bf 4.3. D\' emonstration du th\' eor\`eme 4.2.1.}
        
        \smallskip
        
        On peut supposer que $n$ est pair; en effet, si  $n$  est impair, les invariants de $A_n$ sont d\' etect\' es par $A_{n-1}$ d'apr\`es la prop. 1.6.2.
        
        Supposons donc que $n$ est pair, i.e., $c= 0$ ou $2$, et soient  $ X, E, E(i), X_0, $ comme au §4.1. 
        
        Pour $i = 1,$\dots$,m$, soient $s_i, s'_i, s_i''$ les \' el\' ements non triviaux de $E(i)$.
     Le groupe $E$ est un $2$-groupe ab\' elien \' el\' ementaire de base $s'_1,s_1'',\dots,s'_m,s''_m$. 
        Soit $s_0=1$ si $c=0$, et soit $s_0$ la transposition des deux \' el\' ements de $X_0$ si $c=2$. Soit $s = s_0s_1\!\cdots\!s_{m-1} s_m$; c'est une involution de $S_n$ qui op\`ere sans point fixe sur $X$. Notons $Y$ l'ensemble quotient de $X$ par l'action de $\{1,s\}$; soit $r = |Y| = n/2$. Soit $D$ le centralisateur de $s$ dans $A_n$. Tout \' el\' ement de $D$ d\' efinit par passage au quotient une permutation de $Y$; d'où un homomorphisme $D \to \Sym_Y \simeq S_r$; cet homomorphisme est surjectif, de noyau $E$.
        
        \smallskip
        
        \n {\sc Proposition 4.3.1}. (a) {\it Le groupe $D$ est un groupe de Coxeter
        de type $\sf{D}$$_r$.}
        
        (b) {\it Un \' el\' ement de $D$ est une r\' eflexion si et seulement si c'est une bitransposition de $A_n$ et son image dans $S_r$ est une transposition.}.
        
         (c) {\it L'indice de $D$ dans $A_n$ est impair.}
         
         (d) {\it Le groupe $E$ est un cube maximal de $D$.}
         
      \n    [Rappelons, cf. [Se 22], §4, qu'un {\it cube} d'un groupe de Coxeter fini est un sous-groupe ab\' elien engendr\' e par des r\' eflexions.]
      
      \smallskip
    
\n {\it D\' emonstration}. L'assertion (a) r\' esulte de la construction des groupes de Coxeter de type $\sf{D}$ donn\' ee par exemple dans [Se 23], $\S7$.

 L'assertion (b) se
v\' erifie sur la description des r\' eflexions de $D$ comme r\' eflexions ``longues'' du groupe de type $\sf{B}$ correspondant. 

L'ordre de $D$ est $1\over2$$2^rr!$ et celui de $A_n$
est $1\over2$$n!$. Comme $r = n/2$, on a 

\smallskip

\centerline{$(A_n:D) = n!/2^{n/2}(n/2)! = \prod_{i \ {\rm impair} \ \leqslant n} \ i,$ qui est impair; d'où (c).} 

\smallskip

Le groupe $E$ a pour base  $s'_1,s''_1, $\dots$,s'_m, s''_m$; d'apr\`es (b), ces \' el\' ements sont des r\' eflexions de $D$, ce qui montre que $E$ est un cube de $D$; son rang est $2m$, i.e. $r$ si $r$ est pair, et $r-1$ si $r$ est impair. Dans les deux cas, le rang de $E$ est \' egal au ``rang r\' eduit'' de $D$ au sens de [Se 22], §5.3; d'où (d).

  \medskip
  \n {\it Fin de la d\' emonstration du th\' eor\`eme 4.2.1}.
    
    \smallskip
    
    La partie (c) de la prop. 4.3.1 entraîne que $\Res: \Inv_k(A_n,\cc) \to \Inv_k(D,\cc)$
    est injectif, cf. prop. 1.6.2. D'autre part $D$ est un groupe de Coxeter de type impair (cf. [Se 22], §1.13), ce qui entraîne que ses cubes maximaux sont conjugu\' es ([Se 22], cor. 4.11).
    Tout cube de $D$ est donc contenu dans un conjugu\' e de $E$. D'apr\`es le ``splitting principle'' des groupes de Weyl, cf. [Hi 20] et [GH 22], th. 10 et th. 12, cela entraîne
    que  $\Res: \Inv_k(D, \cc) \to \Inv_k(E,\cc)$ est injectif. D'où l'injectivit\' e de
    $\Res: \Inv_k(A_n, \cc) \to \Inv_k(E,\cc)$.
    
        \medskip
        
        \n {\bf 4.4. Autre exemple de d\' etection lorsque $n \equiv 2$ ou $3$} (mod 4).
        
        \smallskip
        
     \n    {\it Supposons que $n \equiv 2$ ou} $3$ (mod 4).
     
        \smallskip Le plongement $\iota :A_{n-2} \to A_n$ se prolonge en un plongement  $\iota' :S_{n-2} \to A_n$ en d\' efinissant $\iota'(s)$, pour $s \in S_{n-2} \sm A_{n-2}$, comme le produit de $\iota(s)$ et de la transposition $(n-1 \ n)$ de $n-1$ et $n$.  L'image par $\iota'$ d'une
          transposition $(a \  b)$ de $S_{n-2}$ est la
         bitransposition $(a \ b)(n-1 \ n)$ de $A_n$.
         
         Notons $S'_{n-2}$ l'image de $S_{n-2}$ par $\iota'$; c'est un sous-groupe de $A_n$
         d'indice $n(n-1)/2$.
              \smallskip
        
         \n {\sc Proposition 4.4.1}. {\it L'indice de $S'_{n-2}$ dans $A_n$ est impair.}
         
         \smallskip
         
  En effet, $(A_n: S'_{n-2})=n(n-1)/2$, qui est impair car  $n \equiv 2$ ou $3$ (mod 4).
                  
         \medskip
         
         Soit $C$ un cube maximal de $S_{n-2}$, i.e. un $2$-sous-groupe ab\' elien \' el\' ementaire de rang $[(n-2)/2]=2m$, engendr\' e par des transpositions. Soit $C'= \iota'(C)$;  c'est un sous-groupe ab\' elien \' el\' ementaire de $A_n$ de rang $m$, engendr\' e par des bitranspositions.

         \smallskip

         \n {\sc Proposition 4.4.2}.  {\it Les groupes $S'_{n-2}$ et $C'$ d\' etectent les $\cc$-invariants de $A_n$.}

         \smallskip 
         
         L'injectivit\' e de $\Inv_k(A_n,\cc) \to \Inv_k(S'_{n-2},\cc)$ r\' esulte
          de la prop. 4.4.1 et de la prop. 1.5.1. Celle de $\Inv_k(S_{n-2},\cc) \to \Inv_k(C,\cc)$
          r\' esulte de [Se 03], §25 et §29. Celle de $\Inv_k(A_n,\cc) \to \Inv_k(C',\cc)$ r\' esulte des deux pr\' ec\' edentes.

         \smallskip
         
         \n {\it Remarques.}  
         
         1. Soit $C_n$ le cube maximal de $S_n$ engendr\' e par $C$ et la transposition
        $(n-1\ n)$.  On a $C' = C_n \cap A_n$,  ce qui donne une autre caract\' erisation de $C'$.
        
        2. Les groupes $C'$ et $E$ sont isomorphes, mais ne sont pas conjugu\' es dans $A_n$ si $m>0$, i.e. si $n\geqslant 6$. En effet, le nombre de points fixes de $C'$ est $c-2$, alors que celui de $E$ est $c$. Lorsque $n=6$, on peut v\' erifier que $C$ et $E$ sont transform\' es l'un en l'autre par un automorphisme externe de $S_6$; cet automorphisme transforme $q^a$ en $ \<2\>q^a$, cf. [Se 03], (30.11).
                   
\bigskip

    \n {\sc $\S$5. D\' emonstrations des th\' eor\`emes 3.4.1 et 3.4.2}

    \medskip
    
    \n {\bf 5.1. Alg\`ebre \' etales de rang $4$ et $2$-formes de Pfister.}
    
    \medskip
     Soit $K$ un corps de caract\' eristique $\ne 2$. Rappelons qu'une $2$-{\it forme de Pfister} est une forme quadratique de rang 4 sur $K$ qui repr\' esente $1$, et qui est de 
     discriminant 1. Une telle forme peut s'\' ecrire $q= \<1,x,y,xy\> = \<1,x\>\<1,y\>$
     avec $x,y \in K^\times$. Avec les notations de 2.1, on a $\fil(q) \geqslant 2$
     et l'image      de $q$ dans $H^2(K)$ par l'isomorphisme de Milnor est :
     
     \smallskip
    $\ (-x)(-y) = e_k^2+(x)^2+(y)^2+(x)(y) =  e_k^2 + (x)e_k + (y)e_k + (x)(y)$
     
     \smallskip
     \n Un calcul imm\' ediat donne :
        
     \smallskip \n {\sc Lemme 5.1.1}. {\it On a} :
     
     \smallskip
     
  \n    (a) $w_2(q) = (x)^2+(x)(y)+(y)^2.$
     
     \smallskip
\n      (b) $\la^2(q)= 2q-2, \la^3(q)=q, \la^4(q)=1, \la^i(q) = 0$ \ {\it pour} $i >4.$

\medskip
     
     Il est commode de reformuler  (b) en utilisant le polynôme $\la_t(q)= \sum \la^i(q)t^i$, \`a coefficients dans l'anneau $\widehat{W}(K)$ :
     
     \smallskip \n (b$'$) $\la_t(q) = 1 + qt + (2q-2)t^2+qt^3+t^4 = (1+t)^2(1+(q-2)t+t^2).$
     
     \smallskip
     \n {\small [On verra plus loin (th. 6.1.2) que la factorisation de $\la_t(q)$ par une 
     puissance de $1+t$ est un fait g\' en\' eral.]}

\medskip

  Les 2-formes de Pfister sont les formes traces associ\' ees aux $A_4$-torseurs. De façon plus pr\' ecise :
  
  \smallskip
  
  \n {\sc Proposition 5.1.2}. {\it Soit $q$ une forme quadratique sur $K$. Pour qu'il existe 
  un $A_4$-torseur sur $K$ dont la forme trace soit isomorphe \`a $q$,
  il faut et il suffit que $q$ soit une $2$-forme de Pfister.}

\smallskip

\n {\it D\' emonstration.} Si $T$ est un $A_4$-torseur, et si $L$ est l'alg\`ebre \' etale de rang 4 correspondante, la forme $q_L$ repr\' esente 1 puisque $q_L(\frac{1}{2})= 1$,
et son discriminant est 1; c'est donc une 2-forme de Pfister. Inversement, soit
$q =  \<1,x\>\<1,y\>$ une 2-forme de Pfister. Notons $T_x$ le $\Z/2\Z$-torseur
d\' efini par $x$; l'alg\`ebre \' etale correspondante est $L_x=K[X]/(X^2-x)$ et la forme
trace $q_x$ correspondante est $\<2,2x\>$. D\' efinissons de même $T_y, L_y, q_y$.  Le couple $(T_x,T_y)$ d\' efinit un torseur  $T_{x,y}$ du groupe $ \Z/2\Z \times \Z/2\Z$.
L'alg\`ebre \' etale correspondante est l'alg\`ebre biquadratique $L = L_x \otimes L_y$;
sa forme trace est $q_L = q_xq_y = \<2,2x\>\<2,2y\> = \<4\>q = q$. Soit $\sigma$
un isomorphisme entre $ \Z/2\Z \times \Z/2\Z$ et le $2$-Sylow $E$ de $A_4$, cf. §4.1, et soit $T$ le $A_4$-torseur d\' eduit de $T_{x,y}$ par $\sigma$. L'alg\`ebre \' etale correspondante est isomorphe \`a $L$, et sa forme trace est  $q$.

    \medskip
    
    \n {\bf 5.2. Invariants de $A_4$.}  
    
\smallskip 

Lorsque $n=4$, le groupe $E$ du §4.1 est le groupe de type (2,2) dont les trois
    \' el\' ements non triviaux sont les bitranspositions :
      $$s= (1 \ 2)(3 \ 4), \ s' = (1\ 3)(2 \ 4), \ s''= ss' = (1 \ 4)(2\ 3). $$
       On a $A_4 = E\!\cdot\! A_3$ et $S_4=E\!\cdot \!S_3$ (produits semi-directs). 
       
       \smallskip
       
       Identifions $E$ \`a $\Z/2\Z \times \Z/2\Z$ grâce au choix de la base
    $\{s,s'\}$. D'apr\`es [Se 03], §16.4, $\Inv_k(E)$ est un $H(k)$-module libre de rang 4,
   de base $\{1,x,y,xy\}$,   où  $x$ (resp. $y$) est l'\' el\' ement de $\Inv^1_k(E)$ associ\' e \`a la premi\`ere (resp. seconde) projection $E \to \Z/2\Z$, et $xy \in \Inv_k^2(E)$ est le produit de $x$ et de $y$.
      
      \smallskip
   Les groupes $A_3$ et $S_3$ op\`erent sur $E$, donc aussi sur $\Inv_k(E)$.
   
   \smallskip \n    \n {\sc Proposition 5.2.1}. (a) {\it Le sous-$H(k)$-module de $\Inv_k(E)$
   form\' e des \' el\' ements fix\' es par $A_3$ est libre de rang $2$, et a pour base $\{1, ex +ey+xy\}$,   où  $e$ est l'\' el\' ement $e_k=(-1)_k$ de $H^1(k)$, cf.}  §1.2.
   
   (b) {\it  Les \' el\' ements de ce module sont fix\' es par l'action de $S_3$.}
   
   \smallskip {\it D\' emonstration}. L'invariant $ex+ey+xy$ peut aussi s'\' ecrire
   $x^2+y^2+xy$, et ce polynôme quadratique est invariant par $S_3\simeq \GL_2(\F_2)$.    
     Il reste \`a montrer que tout \' el\' ement $z$ de $\Inv_k(E)$ fix\' e par $A_3$ est combinaison
     $H(k)$-lin\' eaire de $1$ et $ex+ey+xy$. Or $z$ s'\' ecrit de façon unique
   sous la forme  $a +bx + cy + dxy$ avec $a,b,c,d \in H(k)$. 
   
   Puisque $z$ est fix\' e par la permutation circulaire $x \mapsto y \mapsto x+y \mapsto x$, on a :
   
   \smallskip
       
  \n \hspace{1mm} $a + bx + cy + dxy = a +by + c(x+y)+ dy(x+y)= a +cx + (b+c+ed)y +dxy$,
  
  \smallskip
   \n  d'où $b= c$ et $c=b+c+ed$, i.e. $b=c =ed$, et $z = a + d(ex + ey + xy)$.
     Cela ach\`eve la d\' emonstration.
     
     \medskip \n {\sc Proposition 5.2.2.} {\it Le $H(k)$-module $\Inv_k(A_4)$ est libre
     de base $\{1, w_2(q^a)\}$.}
     
     \n [Autrement dit, le th. 3.4.1 est vrai pour $n=4$.]
     
     \smallskip \n {\it D\' emonstration.} Comme $E$ d\' etecte les invariants de $A_4$,
     l'homomorphisme de restriction $\Res : \Inv_k(A_4) \to \Inv_k(E)$ est injectif.
     De plus son image est fix\' ee par l'action de $A_3$. Or, d'apr\`es le lemme 5.1.1 (a),
     on a $\Res(w_2(q^a)) = ex +ey + xy$, et, d'apr\`es la prop. 5.2.1 (a), tout \' el\' ement
     de $\Inv_k(E)$ fix\' e par $A_3$  est combinaison $H(k)$-lin\' eaire de $1$ et de 
     $ex+ey+xy$. D'où la proposition.
     
     \smallskip \n {\sc Proposition 5.2.3.} {\it Le $W(k)$-module $\Inv_k(A_4,W)$ est libre de base $\{1, q^a \}$.}
  
     \n [Autrement dit, le th. 3.4.2 est vrai pour $n=4$.]
     
     \smallskip La m\' ethode de d\' emonstration est la même que pour la proposition 
     pr\' ec\' edente. On se ram\`ene \`a d\' eterminer les \' el\' ements de $\Inv_k(E,W)$ qui sont fix\' es par l'action de $A_3$. D'apr\`es [Se 03], th. 27.15, $\Inv_k(E,W)$ a une base form\' ee de quatre \' el\' ements : $1, a_x, a_y, a_{xy}=a_xa_y$  correspondant aux quatre \' el\' ements $1, s, s', ss'$ de $E$. Les trois derniers sont permut\' es circulairement par $A_3$. Il en r\' esulte que le sous-module de $\Inv_k(E,W)$ fix\' e par $A_3$ a pour base $\{1, a_x+a_y+a_{xy}\}$,
donc aussi $\{1, 1 + a_x +a_y + a_{xy}\}$, c'est-à-dire la restriction \`a $E$ de $\{1,q^a\}$.
D'où la proposition.

\medskip

\n {\it Remarque.} Puisque les th\' eor\`emes 3.4.1 et 3.4.2 ont \' et\' e d\' emontr\' es pour $n=4$,
il en est de même des corollaires 3.4.3 et 3.4.4. En particulier, si $\cc = H, W$ ou $ \widehat{W}$, l'homomorphisme $\Res : \Inv_k(S_4,\cc) \to \Inv_k(A_4,\cc)$ est surjectif, et les $\cc$-invariants d'un $A_4$-torseur ne d\' ependent que de la 2-forme de Pfister correspondante.

    \medskip
    
    \n {\bf 5.3. Les $W$-invariants de $A_n$ quand $n$ est multiple de $4$.}
    
\smallskip 
    Supposons que $n = 4m$, avec $m > 0$. Soit $A = \prod_{i=1}^m A(i)$ le sous-groupe de $A_n$ d\' efini au $\S$4.1. Les $A(i)$ sont isomorphes \`a $A_4$; notons $q^a_i$ les formes traces correspondantes, vues comme \' el\' ements de $\Inv_k(A(i),W)$. D'apr\`es la prop. 5.2.3, chaque $\Inv_k(A(i),W)$ est un $W(k)$-module
    libre de base $\{1,q^a_i\}$; de plus, cette base est g\' en\' erique, au sens de 1.8.
    D'apr\`es la prop. 1.8.3 et ses corollaires, on~ a:
    
    \smallskip
   
  \centerline{$\Inv_k(A,W) \ \ = \ \ \Inv_k(A(1),W)\  \otimes \ \Inv_k(A(2),W) \ \otimes \ \cdots \ \otimes \ \Inv_k(A(m),W),$}
  
  \smallskip
\n   où  les produits tensoriels sont pris sur $W(k)$.
  
  \smallskip

  Il en r\' esulte que le $W(k)$-module $\Inv_k(A,W)$ a pour base l'ensemble des
     $q^a_I = \otimes_{i \in I} q^a_i$,    où  $I$ parcourt l'ensemble des parties de
    $\{1,$\dots$,m\}$. Pour tout $d \geqslant 0$, notons $q^a(d)$ la somme des $q^a_I$ tels que    $|I|=d$.
    
\smallskip

    \n {\sc Lemme 5.3.1.} {\it Soit $N$ le normalisateur de $A$ dans $A_n$, et soit $\Inv_k(A,W)^N$ la sous-alg\`ebre de  $ \Inv_k(A,W)$ fix\' ee par $N$. Les $q^a(d)$, $d=0,$\dots$,m$, forment une base de  $\Inv_k(A,W)^N.$}
    
    \smallskip
    
    Cela r\' esulte du fait que l'action de $N$ sur les $m$ facteurs $A(i)$ de $A$ d\' efinit un homomorphisme $N \to S_m$ qui est surjectif, donc permute transitivement les
    $q^a_I$ correspondant \`a des parties $I$ ayant le même nombre d'\' el\' ements.

    \smallskip
    
    Tout \' el\' ement de l'image de $\Res : \Inv_k(A_n,W) \to \Inv_k(A,W)$ est fix\' e par $N$, cf. prop. 1.7.1; d'apr\`es le lemme 5.3.1, c'est une combinaison lin\' eaire des $q^a(d)$.  Ceci s'applique en particulier aux restrictions des \' el\' ements $\la^iq^a$ de $\Inv_k(A_n,W)$.
    Pour \' enoncer le r\' esultat, nous utiliserons (comme au $\S$5.1) le polynôme $\la_t(q^a)$ en une variable $t$, d\' efini par la formule $\la_t(q^a)=\sum_{i \geqslant 0} t^i\la^iq^a$; c'est un polynôme \`a coefficients dans l'anneau $\Inv_k(A,W)$.
     
     \smallskip
     
     \n{\sc Proposition 5.3.2}. {\it On a  $\la_t(\Res(q^a)) = (1+t)^{2m} \sum_{d \leqslant m} t^d(1-t)^{2m-2d}q^a(d)$ dans l'anneau de polynômes} $\Inv_k(A,W)[t]$. 
    
       \smallskip
   
   \n {\it D\' emonstration.} 
   
   \smallskip \n On a $\Res(q^a) = \sum_i q^a_i$, d' où  $\la_t(\Res(q^a)) =  \prod_i \la_t(q^a_i).$
   
   \medskip
   \n Comme  $\la_t(q^a_i) = (1+t)^2(1+(q^a_i-2)t+t^2)$ (cf §4.1), on obtient:
   
   \smallskip
 \centerline{ $\la_t(\Res(q^a)) = (1+t)^{2m}\prod_i (1+(q^a_i-2)t+t^2).$}
   
         \medskip \n  Le facteur $\prod_i (1+(q^a_i-2)t+t^2)$ est \' egal à
$\prod_i (\alpha + \beta q^a_i)$,   où  $\alpha = (1-t)^2, \beta = t$. 

\smallskip
\n Vu la d\' efinition des $q^a(d)$, il est \' egal \`a 

\smallskip
$\sum_{d\leqslant m}\alpha^{m-d}\beta^dq^a(d) = \sum_{d\leqslant m}(1-t)^{2m-2d}t^dq^a(d)$,

\smallskip
\n d'où la proposition.

        \smallskip

     \n {\sc Corollaire 5.3.3.} {\it Soit $d$ tel que $0\leqslant d \leqslant m$. On a} :
     
     \smallskip
  \centerline{ $\Res(\la^dq^a) = q^a(d) \ + $ {\it combinaison $\Z$-lin\' eaire des} $q^a(d'), d' < d.$}
     
     \smallskip
         
        \n   En effet, le coefficient de $q^a(d)$ dans le polynôme de la prop. 5.3.2 est un polynôme en $t$ dont tous les termes sont de degr\' e $>d$, \`a la seule exception de~$t^d$.

          \smallskip
        
        \n {\sc Proposition 5.3.4.} {\it Les invariants $\Res(\la^dq^a), 0 \leqslant d \leqslant m$, forment une $W(k)$-base de $\Inv_k(A,W)^N.$}

          \smallskip
        \n En effet, d'apr\`es le cor. 5.3.3, ils se d\' eduisent de la base des $q^a(d)$ par une  matrice triangulaire avec des 1 sur la diagonale.
        
        \smallskip
        
        \n {\sc Corollaire 5.3.5.} {\it L'homomorphisme \ $\Res: \Inv_k(A_n,W) \to \Inv_k(A,W)$ d\' efinit un isomorphisme de $\Inv_k(A_n,W)$ sur $\Inv_k(A,W)^N$.}
        
        \smallskip
        
     \n En effet, cet homomorphisme est injectif (cor. 4.2.2), et il est surjectif puisque
     son image contient une base de $\Inv_k(A,W)^N$ d'apr\`es la prop. 5.3.4.
              
            \smallskip
        
        \n {\sc Corollaire 5.3.6.} {\it Le $W(k)$-module $\Inv_k(A_{4m},W)$ est un module libre de base $\{\la^0(q^a),$\dots$, \la^m(q^a)\}$.}
        
        \n [Autrement dit, le th. 3.4.2 est vrai lorsque $n=4m$.]
        
        \smallskip
        \n Cela r\' esulte du cor. 5.3.5 et de la prop. 5.3.4.

        \medskip

    \n {\bf 5.4. Les $H$-invariants de $A_n$ quand $n$ est multiple de $4$.}
    
\smallskip 
       La m\' ethode pour les $H$-invariants est la même que celle du §5.3 pour les
       $W$-invariants. On utilise encore les groupes $A$ et $A(i)$. Si  $I$ est une partie
       de $\{1,m\}$, posons $w(I) = \otimes_{i\in I} w_2(q^a_i)$; c'est un \' el\' ement de
       $\Inv^{2|I|}_k(A)$; d'apr\`es le cor. 1.8.2 et la prop. 1.8.3, les $w(I)$ forment une $H(k)$-base de $\Inv_k(A)$. Pour tout $d \geqslant 0$, soit $w(d)$ la somme des $w(I)$
       avec $|I| = d$. L'analogue du lemme 5.3.1 est l'\' enonc\' e suivant, qui se d\' emontre de
       la même mani\`ere :

       \smallskip 
    \n {\sc Lemme 5.4.1.} {\it Soit $N$ le normalisateur de $A$ dans $A_n$, et soit $\Inv_k(A)^N$ la sous-alg\`ebre de  $ \Inv_k(A)$ fix\' ee par $N$. Les $w(d)$, $d=0,$\dots$,m$, forment une base de  $\Inv_k(A)^N.$}
    
    \smallskip
    
   La prop. 5.3.2 et le cor. 5.3.3 sont remplac\' es par l'\' enonc\' e suivant :

       \smallskip 
    \n {\sc Lemme 5.4.2.} {\it On a $\Res(w_{2d}(q^a)) = w(d)$ pour tout $d\geqslant 0$.}
   
   \smallskip
    Cela r\' esulte du fait que $\Res(q^a)= \Sigma \  q^a_i$ et que la classe de Stiefel-Whitney totale de $q^a_i$ est $1 + w_2(q^a_i)$.
    
    \smallskip
    
      On d\' eduit de là, par les mêmes arguments qu'au §5.3:

          \smallskip
        
        \n {\sc Proposition 5.4.3.} {\it Les invariants $\Res(w_{2d}(q^a)), 0 \leqslant d \leqslant m$, forment une $H(k)$-base de $\Inv_k(A)^N.$}

        \smallskip
        
        \n {\sc Corollaire 5.4.4.} {\it L'homomorphisme \ $\Res : \Inv_k(A_n) \to \Inv_k(A)$ d\' efinit un isomorphisme de $\Inv_k(A_n)$ sur $\Inv_k(A)^N$.}
        
        \smallskip
        
        \n {\sc Corollaire 5.4.5.} {\it Le $H(k)$-module $\Inv_k(A_{4m})$ est un module libre de base les $w_{2d}(q^a), \ 0 \leqslant d \leqslant m$.}
        
        \n [Autrement dit, le th. 3.4.1 est vrai lorsque $n=4m$.]

    \medskip
      
    \n {\bf 5.5. Fin de la d\' emonstration des th\' eor\`emes $3.4.1$ et $3.4.2$.}
    
    \smallskip

\n {\it Le cas du th\' eor\`eme 3.4.1.}

  On raisonne par r\' ecurrence sur $c= n - 4m$,   où  $m=[\frac{n}{4}]$. On a
 $c = 0, 1, 2$ ou $3$. Le cas $c=0$ r\' esulte du cor. 5.4.5. Supposons $c>0$.
 D'apr\`es l'hypoth\`ese de r\' ecurrence, $\Inv_k(A_{n-1})$ a pour base les $w_{2d}(q_{n-1}^a)$ pour $d = 0,$\dots$,m$. L'homomorphisme de restriction $\Res : \Inv_k(A_n)\to \Inv_k(A_{n-1})$ transforme $q_n^a$ en $q_{n-1}^a+1$. D'où
$\Res(w_{2d}(q_n^a))= w_{2d}(q_{n-1}^a)$ pour tout $d$. Cela montre que $\Res : \Inv_k(A_n) \to \Inv_k(A_{n-1}$ est surjectif, donc bijectif d'apr\`es le cor. 4.2.3. On en conclut que les $w_{2d}(q_n^a)$ forment une base de $\Inv_k(A_n)$.

\smallskip

\n {\it Le cas du th\' eor\`eme 3.4.2.}

La m\' ethode est la même : le cas $c=0$ r\' esulte du cor. 5.3.6, et l'on raisonne par r\' ecurrence sur $c$ quand $c>0$. La seule diff\' erence consiste en la formule donnant $\Res(\la^d(q^a_n))$,
qui est
   $$\Res(\la^d(q^a_n)) =\la^d(q^a_{n-1}) + \la^{d-1}(q^a_{n-1}),$$
   
\n   puisque $\Res(q_n^a) = q_{n-1}^a + 1$. 

D'apr\`es l'hypoth\`ese de r\' ecurrence, les
$\la^d(q^a_{n-1}) + \la^{d-1}(q^a_{n-1})$, $d=0,$\dots$,m$, forment une base de $\Inv_k(A_{n-1},W)$. 

On en conclut que $\Res : \Inv_k(A_n,W) \to \Inv_k(A_{n-1},W)$ est bijectif, et que les $\la^d(q^a_n)$ forment une base de
$\Inv_k(A_n,W)$. 

\medskip

Ceci ach\`eve la d\' emonstration des \' enonc\' es du §3.4. 

\medskip

Explicitons une cons\' equence des th\' eor\`emes 3.4.1, 3.4.2 et 4.2.1. Pour l'\' enoncer, disons qu'une alg\`ebre \' etale $L$ de degr\' e $n$ sur un corps $K$ {\it est de type} (T)  (initiale de ``test'') si c'est un produit de $c$ copies de $k$ et de $m$ alg\`ebres de degr\' e 4 qui  sont {\it biquadratiques}, i.e. produit tensoriel de deux alg\`ebres \' etales de degr\' e 2.
Cela entraîne que $q_L= c + \sum_{i=1}^m q_i$,  où  les $q_i$
sont des $2$-formes de Pfister. Le fait que $E$ d\' etecte les invariants de $A_n$ équivaut à :

\smallskip
    
        \n {\sc Th\' eor\`eme 5.5.1.} {\it Soit $a \in \Inv_k(A_n,\cc)$ tel que $a(L)=0$
        pour toute alg\`ebre \' etale $L$ de degr\' e $n$ de type $(T)$ sur toute extension $K$ de $k$. Alors $a=0$.}
                 
       \bigskip

    \n {\sc $\S$6. Propri\' et\' es des formes traces associ\' ees aux $A_n$-torseurs}

    \medskip

       \n {\bf 6.1. Relations entre les invariants $\la^iq^a$.}
       
       \smallskip
       
       Conservons les notations $k, n, m, c$ du $\S$5.5. Soit $L$ une
       $k$-alg\`ebre \' etale de dimension $n$ et de discriminant 1, et soit $q_L$ sa forme 
       trace. D'apr\`es le th. 3.4.2, les $\la^jq_L$ sont des combinaisons
      lin\' eaires de ceux avec $j$ est $\leqslant m$. On va pr\' eciser
      cet \' enonc\' e :

          \smallskip
        
        \n {\sc Th\' eor\`eme 6.1.1.} {\it Il existe des entiers $z(i,j,n)$ tels que
         
         \medskip
        \hspace{20mm} $\la^jq_L = \sum_{i\leqslant m}z(i,j,n)\la^ iq_L$ \ dans \ $\widehat{W}(k)$
         
         \medskip
        \n   pour tout couple $(k,L)$ comme ci-dessus.}
             
        \smallskip
      
     La d\' emonstration du th. 6.1.1 sera donn\' ee au §6.2. Elle utilise le th\' eor\`eme suivant :
     
     \smallskip
       \n {\sc Th\' eor\`eme 6.1.2.} {\it Soient $k, n, m,c, L $ comme ci-dessus. Le polynôme
       $\lambda_t(q_L)$ est divisible par $(1+t)^{2m+c}$}.
       
       \n  [Il s'agit de polynômes \`a coefficients dans $\widehat{W}(k).$]
       
       \smallskip
       
      \n {\it D\' emonstration du th. 6.1.2.}  Remarquons d'abord que le th\' eor\`eme est vrai
      si $L$ est de type (T) au sens du §5.5. En effet, $q_L$ est alors somme de
      $2$-formes de Pfister $q_i$\ $ (i=1,\dots,m)$, et de la forme $\<1, 1,\dots,1\>$ de rang $c$, et l'on a 
      $$\la_t(q_L) = (1+t)^c\prod_i \la_t(q_i) = (1+t)^{2m+c}\prod_i(1 + (q_i-2)t+t^2).$$
        \n  Nous allons nous ramener \`a ce cas.
      
     Notons d'abord qu'un polynôme $f(t)= \sum _ia_it^i$ \`a coefficients dans un anneau commutatif est divisible par une  puissance $(1+t)^N$ de $1+t$ si et seulement si {\it ses d\' eriv\' ees divis\' ees $\Delta^jf(t) = \sum_{i\geqslant j} {i\choose j} a_jt^{i-j}$ s'annulent pour $t=-1$ quel que soit  $j<N$}: c'est une cons\' equence de la formule de Taylor pour les polynômes, cf. [A IV], §IV.8. Cela montre que le th. 6.1.2 est \' equivalent au suivant :
           
      \smallskip
       \n {\sc Th\' eor\`eme 6.1.3.} {\it  Pour tout $j < 2m+c$, on a} $ \sum_{i\geqslant j} (-
1)^i{i\choose j}\la^{i}(q_L) = 0.$

\medskip   
      \n {\it Exemples} (où l'on \' ecrit $\la^i$ \`a la place de $\la^i(q_L)$).
      
       \smallskip
      
  Le cas $j=0$ donne : \    $\la^0 -\la^1 + \la^2- \la^3 + \la^4 - \dots = 0$  \quad  si $n \geqslant 1$.
      
       \smallskip
  Le cas $j=1$ donne : \   $ \la^1 - 2\la^2 +3\la^3 -4 \la^4 \ + \ \dots   = 0$ \ \ \ si $n \geqslant 2.$

       \smallskip
       
       \n {\it D\' emonstration du th. 6.1.3.} 
       
       Soit $j < 2m+c$. Soit 
      $ a_j = \sum_{i\geqslant j} (-1)^i{i\choose j}\la^i(q^a)$; c'est un $\widehat{W}$-invariant de $A_n$ qui s'annule pour toute alg\`ebre \' etale de type (T).
      D'apr\`es le th. 5.22, on a $a_j=0$. D'où le th\' eor\`eme. 
      
      \smallskip
      \n Cela ach\`eve la d\' emonstration du th. 6.1.3, et donc aussi du th. 6.1.2.
      
      \medskip \n {\it Remarque}. On peut aussi d\' eduire le th. 6.1.2 du th. 29.4 de [Se 03].
                        
       \medskip 
       
       \n {\bf 6.2. D\' emonstration du th\' eor\`eme 6.1.1.}
            
            \smallskip
            
 D'apr\`es le th.6.1.2, on a $\la_t(q_L) = (1+t)^{2m+c}P(t)$,  où  $P(t)$ est un polynôme de degr\' e $2m$ \`a coefficients dans $\widehat{W}(k)$. Soient $p_0,...,p_{2m}$ les coefficients de $P$.
Les $\la^i(q_L)$ sont des combinaisons $\Z$-lin\' eaires des $p_j$ avec $j\leqslant i$ :
   
   \medskip 
\centerline{(*) \quad $\la^i(q_L) = \sum_{j \leqslant i}$$2m+c \choose i-j$$ p_j.$}

\smallskip
    
    Inversement, la formule $P(t) =(1+t)^{-2m-c} \la_t(q_L)$, valable dans l'anneau de s\' eries formelles  $\widehat{W}(k)[[t]]$, montre que les $p_i$ sont combinaisons $\Z$-lin\' eaires des $\la^j(q_L)$ avec $j\leqslant i$:
    
    \smallskip    \centerline{$p_i = \sum_{j \leqslant i}$$-2m-c \choose i-j$$ \la^j(q_L).$}
    
    \smallskip
    En particulier, les $p_i$ tels que $i\leqslant m$ sont des combinaisons $\Z$-lin\' eaires explicites de $\la^0(q_L), \dots, \la^m(q_L)$. Or le polynôme $\la_t(q_L)$
    est un polynôme {\it r\' eciproque} : on a $\la^i = \la^{n-i}$ pour tout $i$. Il en est donc de même de $P(t)$: on a $p_i = p_{2m-i}$ pour tout $i$. En appliquant ceci 
    aux $i > 2m$, on voit que tous les $p_i$ sont des combinaisons $\Z$-lin\' eaires explicites de $\la^0(q_L), \dots, \la^m(q_L)$. D'apr\`es (*), il est de  même
    de tous les $\la^i(q_L)$. Cela d\' emontre le th. 6.1.1. 

\medskip

\n {\it Exemple}. Traitons le cas $n=7$, en \' ecrivant $\lambda^i$ \`a la place de $\lambda^i(q_L)$. On a $m=1, c=3$ de sorte que
$\la^2, \la^3, \dots, \la^7$ sont des combinaisons $\Z$-lin\' eaires de $\la^0 = 1$ et $\la^1= q_L$. Plus pr\' ecis\' ement :
$$ \la^2= \la^5= 5\la^1-14,  \quad \la^3=\la^4 = 10\la^1-35, \quad \la^6 =\la^1, \quad \la^7=1.$$
\n  En effet, on a $\la_t(x)=(1+t)^5P(t)$,  où  
 $P(t) =1+xt+t^2$, avec $x \in \widehat{W}(k)$, autrement dit:
 $$\sum_{i=0}^7 \la^i t^i = (1+5t+10t^2+10t^3+5t^4+t^5)(1+xt+t^2).$$
\n En comparant les termes en $t$, on obtient $x=\la^1-5$; les termes en $t^2$ et $t^3$  donnent
$\la^2$ et $\la^3$, et l'on en d\' eduit les autres $\la^i$ puisqu'ils sont \' egaux aux $\la^{7-i}$.

\newpage
 \n{\bf Appendice. Une g\' en\' eralisation de la formule $(1.1.5)$ du $\S$1.1.} 

\bigskip
\n {\bf A.1. Hypoth\`eses.}

\smallskip
Soit $R = \oplus_{d\geqslant 0} R_d$ un anneau commutatif gradu\' e. Supposons
qu'il ait la propri\' et\' e suivante:

\smallskip
(A.1.1) {\it Pour tout $d\geqslant 0$, et tous $x,y \in R_d$, on a $x^2y+xy^2=0$.}

\smallskip Noter que, pour $(d,x,y)=(0,1,1)$, cette condition entraîne que $2=0$ dans $R$, autrement dit que $R$  est une $\F_2$-alg\`ebre. Il en r\' esulte que la 
propri\' et\' e $x^2y+xy^2=0$ de (A.1.1) peut se r\' ecrire $x^2y=xy^2$, ou  aussi :

\smallskip

(A.1.2) \ $(1+x)(1+y) = (1+x+y)(1+xy)$.

\smallskip

\smallskip \n {\it Exemples.}  1) Soit $H$ une $\F_2$-alg\`ebre gradu\' ee commutative contenant un \' el\' ement $\e$ de degr\' e 1 tel que $x^2= \e^{d}x$
pour  tout $d \geqslant 0$ et tout $x \in H_d$. Alors $H$ a la propri\' et\' e (A.1.1). Ceci s'applique en particulier \`a l'alg\`ebre de cohomologie $H(k)$ du $\S$ 1.1.

2) Si $V$ est un $\F_2$-espace vectoriel, l'alg\`ebre ext\' erieure $\wedge V$
 a la propri\' et\' e (A.1.1).
 
 3) L'alg\`ebre de polynômes $\F_2[X]$ a la propri\' et\' e (A.1.1) - c'est le cas 
 particulier $k = \R$ de l'exemple 1).

\medskip
\n {\bf A.2. Enonc\' e du th\' eor\`eme A.2.3.}

\medskip
\n Pour tout $n\geqslant 1$, notons $P_n$ l'ensemble des \' el\' ements de $R$ de la forme

\medskip
 \ (A.2.1)\quad \ \ \ $x = \prod_{i=0}^{n-1} (1+ a_i)$, \ avec  $a_i \in R_{2^i}$
  pour tout $i$.

\medskip
\n On a $P_1 = 1 + R_1, P_2 =P_1\!\cdot\!(1+R_2)$, et $P_{n+1} = P_n\!\cdot\!(1+ R_{2^n})$ pour tout $n$.

\medskip

{\it Remarque.} Soit $x$ comme dans (A.2.1), et soit  $n \geqslant 0$, \' ecrit sous forme dyadique $n = \sum_{i\in I} 2^ i$, où $I$ est une partie finie de $\N$. La $n$-i\`eme composante $x_n$ de $x$ est donn\' ee par la formule :

\smallskip
 (A.2.2)  $x_n = \prod_{i\in I} a_i$.

\smallskip \n En particulier, on a $x_{2^i} = a_i.$

\medskip
     \n {\sc Th\' eor\`eme A.2.3.} {\it On a $P_n \! \cdot\!  P_n \subset P_{n+1}$.}

\smallskip
Les $P_n$ forment une suite croissante de parties de $R$; notons $P$
leur r\' eunion. Le th\' eor\`eme A.2.3. implique:

\smallskip
\n {\sc Corollaire A.2.4.} {\it On a $P \! \cdot\!  P \subset P$.}

Autrement dit, $P$ est stable par multiplication ; noter que ce n'est pas en g\' en\' eral un groupe, car ses \' el\' ements ne sont pas toujours inversibles dans $R$. Par exemple, si $R= \F_2[X]$, aucun \' el\' ement de $P$, \`a part 1, n'est inversible.

\medskip
\n {\bf A.3. D\' emonstration du th\' eor\`eme A.2.3.}

\smallskip
\n {\sc Lemme A.3.1.} {\it Soit $x \in P_n$ et soit $y=1+z$, avec $z \in R_{2^m}, m \geqslant 0$.  On a $xz \in P$.}

\smallskip

\n {\it D\' emonstration}. Ecrivons $x$ sous la forme (A.2.1) :

\smallskip

 \n \quad \ \ \ $x = \prod_{i=0}^{n-1} (1+ a_i)$, \ avec $\deg(a_i) = 2^i$ pour tout $i$.

\smallskip

Si $m \geqslant n$, le produit $xy$ appartient \`a $P_{m+1}$, donc \`a $P$.
On peut donc supposer dans le lemme A.3.1 que $m\leqslant n$.
Nous proc\`ederons par r\' ecurrence sur $ n-m$, autrement dit nous
 supposerons que le lemme est vrai pour les quadruplets $(n',m',x',z')$ avec $n'-m' < n-m$.

 \smallskip
  D\' ecomposons  $x$  en  $x = x_-x_o x_+$,  où  :

\smallskip
\centerline {$x_- = \prod_{i<m}(1+a_i), \ x_o = 1+ a_m, \ x_+ = \prod_{i>m}(1+a_i).$}

\smallskip
On a $xy = x_-(1+a_m)(1+z)x_+ = x_-(1+a_m+z)(1+a_mz)x_+$,
autrement dit $xy = uv$ avec $u=x_-(1+a_m+z)$ et $v =(1+a_mz)x_+.$

\smallskip
Le produit  $u=x_-(1+a_m+z)$ appartient \`a $P_{m+1}$. D'autre  part, l'hypoth\`ese de r\' ecurrence, appliqu\' ee au quadruplet $(n,m+1,x_+,a_mz)$,
montre que  $v$ appartient \`a $P$; de plus ses composantes
de degr\' e $\delta$ tel que $0 <  \delta < 2^{m+1}$ sont 0. Autrement dit,
$v$ est produit de facteurs du type $(1+b_j)$ avec $\deg(b_j) = 2^j, \ j  \geqslant m+1.$  Le produit $uv$ appartient donc \`a $P$.

\smallskip

\n {\it Fin de la d\' emonstration du th\' eor\`eme A.2.3.}

Le lemme A.3.1 montre que $P$ est stable par multiplication par les \' el\' ements de la forme $1+z$, avec $\deg(z)$ puissance de 2. Comme ces \' el\' ements engendrent $P$, on a donc $P \! \cdot\!  P \subset P$. En particulier, on a $P_n \! \cdot\!  P_n \subset P$. De plus les
composantes des \' el\' ements de $P_n$ sont nulles en tout degr\' e $> 2^n$;
celles des \' el\' ements de  $P_n \! \cdot\!  P_n$ sont donc nulles en degr\' e
$> 2^{n+1}$; cela montre que  $P_n \! \cdot\!  P_n$ est contenu dans $P_{n+1}$.

\smallskip

\n {\it Exemple: le cas $n= 2$.} 

Soient $(1+a_0)(1+a_1)$ et $(1+b_0)(1+b_1)$ deux \' el\' ements de $P_2$ (avec $a_0, b_0 \in R_1$, et $a_1,b_1 \in R_2$). Leur produit est
l'\' el\' ement $ (1+c_0)(1+c_1)(1+c_2)$ de $P_3$,  où  :

\smallskip

 $c_0 = a_0+b_0$,
 
 $c_1 = a_1+ a_0b_0 + b_1$,
   
       $c_2 = a_0a_1b_0 + a_1b_1+a_0b_0b_1.$
   
   \medskip    
   
\n {\bf A.4. D\' emonstration de la formule (1.1.5) du $\S$ 1.1.}

 \smallskip
 
 Rappelons de quoi il s'agit.

\smallskip
Soit $d$ un entier $\geqslant 1$, et soient $y_i$ des \' el\' ements de $R_d$,
en nombre fini. Pour tout  $m \geqslant 1$, soit $s_m$ la $m$-i\`eme fonction sym\' etrique \' el\' ementaire des $y_i$; on a
$$  \ \ ({\rm A}.4.1) \quad 1 + \sum s_m = \prod(1+y_i). \hspace{7cm} $$
La formule (1.1.5) affirme que
$$  ({\rm A}.4.2) \quad 1 + \sum s_m = \prod_{j\geqslant 0}(1 + s_{2^j}). \hspace{6.5cm} $$

\smallskip
\n {\it D\' emonstration de} (A.4.2) {\it lorsque $d=1$.}

\smallskip  On a $1 + y_i \in P_1$ pour tout $i$. D'apr\`es le th\' eor\`eme A.2.3, cela entraîne que $\prod (1+y_i)$ appartient \`a $P$, donc peut s'\' ecrire 

\smallskip  
  (A.4.3) \quad $ \prod(1+y_i) = \prod (1+ a_j)$ avec $a_j \in R_{2^j}$ pour tout $j \geqslant 0$.

\smallskip
Cette formule montre que $a_j$ est la fonction sym\' etrique \' el\' ementaire des
  $y_i$ d'indice $2^j$, autrement dit est \' egale \`a $s_{2^j}$. 
   
   \medskip

\n {\it D\' emonstration de} (A.4.2) {\it dans le cas g\' en\' eral.}

\smallskip

Soit $R[d] = \oplus_{n \geqslant 0} R_{dn}$; c'est un sous-anneau de $R$.
Munissons-la de la graduation  où  les \' el\' ements de $R_{dn}$ sont de degr\' e
$n$, autrement dit $R[d]_n = R_{dn}$.  Cet anneau v\' erifie lui aussi la condition $x^2y = xy^2$ lorsque $x,y$
ont le même degr\' e. Les $y_i$ sont des \' el\' ements de degr\' e 1 de $R[d]$. On peut donc leur appliquer (A.4.2), et cela donne le r\' esultat cherch\' e.

   \medskip
\n {\bf A.5. Produits infinis et s\' eries formelles.}   
    
      \smallskip
Dans la th\' eorie des classes de Stiefel-Whitney, il est souvent commode
d'introduire des s\' eries formelles, not\' ees $a_t$, qui sont des \' el\' ements 
de $R[[t]]$ : 

    \smallskip
 $a_t = \sum_{n \geqslant 0} a_nt^n$ telles que $a_n \in R_n$ et $a_0=1$.

    \smallskip
Ces s\' eries forment un groupe $\Pi$ pour la multiplication; ce groupe est complet pour la topologie naturelle de $R[[t]]$ (où $R$ est muni de la topologie discr\`ete).

    \smallskip

\n {\bf Th\' eor\`eme A.5.1}. {\it Les \' el\' ements de $\Pi$ de la forme}

  \smallskip
(A.5.2) \ \ $a_t = \prod_{j \geqslant 0} (1+ t^{2^j}a_j)$, {\rm} avec $a_j \in R_{2^j}$,

  \smallskip
\n{\it forment un sous-groupe ferm\' e $\Pi^1$ de $\Pi$}.

       \smallskip
   
   \n {\it D\' emonstration.}  Si $a_t \in \Pi^1$ est un polynôme en $t$, alors 
   $\prod_j(1+ a_j)$ est un \' el\' ement de $P$. On a donc un plongement naturel de $P$ dans $\Pi^1$, compatible au produit. Il est clair que $P$ est dense dans $\Pi^1$. D'apr\`es le th\' eor\`eme A.2.3, $P$
   est stable par la multiplication. Il en est donc de même de $\Pi^1$.
   De plus, les \' el\' ements de $P$ sont inversibles dans $\Pi^1$ : cela r\' esulte 
 de la formule:
 
   \smallskip
 (A.5.2) \ \ $(1+x)^{-1} = (1+x)(1+x^2)(1+x^4)\cdots $.   
   
     \smallskip \n Cela entraîne que $\Pi^1$ est un sous-groupe de $\Pi$.
   
     \smallskip
  \n {\it Remarque}. Pour tout $d > 0$ on peut d\' efinir de façon analogue
  un sous-groupe ferm\' e $\Pi^d$ de $\Pi$; le groupe $\Pi$ est produit direct
  des $\Pi^d$, pour $d$ impair.

          \centerline{ {\bf R\' ef\' erences}}
\vspace{4mm}
 
          \n [AP 71] J. Arason \& A. Pfister, {\it Beweis des Krullschen Durschnitzzatses f\' eür
     den Wittring}, Invent. math. {\bf12} (1971), 173--176.
     
       \n  [A IV-VII] N. Bourbaki, {\it Alg\`ebre. Chapitres} IV-VII, Masson, Paris, 1981.

    \n [BS 94] E. Bayer-Fluckiger \& J-P. Serre, {\it Torsions quadratiques et bases normales autoduales}, Amer. J. Math. {\bf116} (1994), 1--64.
        
     \n [BS 21] -----, {\it  Lines on cubic surfaces, Witt invariants and Stiefel-Whitney classes}, Indag. Math. (N.S.) {\bf32} (2021), 920--938.

     \n [Di 02] L.E. Dickson, {\it Theorems on the residues of multinomial coefficients with respect to a prime modulus}, Quarterly
     Journal Pure Applied Mathematics {\bf33} (1902), 378--384.
     
\n [Ga 20] N. Garrel, {\it Witt and cohomological invariants of Witt classes}, Ann. K-Theory {\bf 5} (2020), 213-248.

\n [Ga 24] -----, {\it Even Stiefel-Whitney invariants for anti-hermitian quaternionic forms}, arXiv 2304.07807.
 
\n [GH 22] S. Gille \& C. Hirsch, {\it On the splitting principle for cohomological invariants of reflection groups}, Transformation Groups {\bf27} (2022), 1261--1285.

\n [GKT 89] J. Gunarwardena, B. Kahn \& C. Thomas, {\it Stiefel-Whitney classes of real 
representations of finite groups}, J. Algebra 126 (1989), 327--347.
  
 \n [Hi 09] C. Hirsch, {\it Cohomological invariants of reflection groups}, Diplomarbeit (Betreuer: Prof. Dr. Fabien Morel), Univ. München, 2009; version révisée : arXiv 1805.04670.
  
\n [Hi 20] -----, {\it On the decomposability of mod $2$ cohomological invariants
of Weyl groups}, Comment. Math. Helv.  {\bf95} (2020), 765--809.

\n [Ka 84] B. Kahn, {\it Classes de Stiefel-Whitney de formes quadratiques et de repr\' esentations galoisiennes r\' eelles}, Invent. math. 78 (1984), 223-256.

\n [Ka 20]  -----, {\it Comparison of some field invariants}, J. Alg. {\bf232} (2000), 485--492.

  \n [KMRT 98]  M-A. Knus, A. Merkurjev, M. Rost \& J-P. Tignol, {\it The Book of Involutions}, A.M.S. Colloq. Publ. {\bf44}, 1998.
  
\n [La 05] T.Y. Lam, {\it Introduction to quadratic forms over fields}, A.M.S. GSM {\bf67}, 2005.

\n [Me 10] A. Merkurjev, {\it Developments in algebraic K-theory and quadratic forms after the work of Milnor}, in J. Milnor, {\it Collected Papers} V , 399--418.

\n [Mi 70]  J. Milnor, {\it Algebraic K-theory and quadratic forms}, Invent. math. {\bf9} (1969/1970), 318-344; Coll. Papers V, 347--373.

\n [Mo 05] F. Morel, {\it  Milnor's conjecture on quadratic forms and mod $2$ motivic complexes}, Rend. Sem. Mat. Univ. Padova {\bf114} (2005), 63--101.

\n [OVV 07]  D. Orlov, A. Vishik \& V. Voevodsky, {\it An exact sequence for $KM_{*}/2$ with applications to quadratic forms}, Ann. of Math. (2) {\bf165} (2007), 1--13.

\n [Se 65] J-P. Serre, {\it Cohomologie Galoisienne}, LNM 5, Springer-Verlag, 1965; cinqui\`eme \' edition r\' evis\' ee et compl\' et\' ee, 1994; trad. anglaise {\it Galois Cohomology}, Springer-Verlag, 1997.

\n [Se 03] -----, {\it Cohomological invariants, Witt invariants, and trace forms}, notes by Skip Garibaldi, ULS {\bf28}, A.M.S, 2003, 1--100; version révisée dans Coll. Papers V, n° 179.

\n [Se 18] -----, {\it Cohomological invariants {\rm mod} $2$ of Weyl groups}, Oberwolfach Report {\bf21} (2018), 1284--1286; arXiv :1805.07172; Coll. Papers V, n° 204.

\n [Se 22] -----, {\it Groupes de Coxeter finis $:$ involutions et cubes}, L'Ens. Math. {\bf68} (2022), 99--133, Coll. Papers V, n° 209.

\n [SGA 6]  P. Berthelot, A. Grothendieck  \& L. Illusie, {\it Th\' eorie des intersections et th\' eor\`eme de Riemann-Roch},  S\' eminaire de G\' eom\' etrie Alg\' ebrique du Bois-Marie, LNM {\bf225} (1971).

\n [Sp 52] T.A. Springer, {\it Sur les formes quadratiques d'indice z\' ero}, C.R.A.S. {\bf234} (1952), 1517--1519.

\n [To 22] B. Totaro, {\it Divided powers in the Witt ring of symmetric bilinear forms}, Ann.  K. Theory {\bf8} (2023), 275--284.

\n [Vi 09] C. Vial, {\it Operations in Milnor theory}, J. Pure Applied Algebra {\bf13} (2009),  1325--1345.

\n [Vo 03.I] V. Voevodsky, {\it Reduced powers in motivic cohomology}, Publ. Math. I.H.E.S. {\bf98} (2003), 1-57.

\n [Vo 03.II]  -----, {\it Motivic cohomology with $\Z/2$ coefficients}, Publ. Math. I.H.E.S. {\bf98} (2003), 59--104.

\bigskip

\n Coll\`ege de France, 3, rue d'Ulm, 75005 Paris

       \end{document}